\documentclass{article}
\usepackage[T1]{fontenc}
\usepackage[utf8]{inputenc}

\usepackage{graphicx}
\usepackage{fullpage, amsmath, amssymb, amsthm}
\usepackage{caption}
\usepackage{subcaption}

\usepackage{url}
\usepackage{graphicx, physics}
\usepackage{subcaption}

\title{Explaining the poor performance of the KASS algorithm implementation}
\author{Allan Grønlund\\
Computer Science, Aarhus University, jallan@cs.au.dk}
\date{Draft: First Version December 13th 2019, Second Version January 28, 2020, \\Current Version March 12th, 2020}
\begin{document}
\maketitle
\begin{abstract}
  By investigating the code for the KASS algorithm implementation used in the paper ``Exploring the quantum speed limit with computer games'' \cite{sorensen2016} by S{\o}rensen et al. (provided by the authors),
  we describe how the poor performance of the KASS algorithm reported in \cite{sorensen2016} is entirely caused by a simple sign error in a derivative calculation.
  Changing only this one sign in the KASS implementation, we show that the algorithm provides results comparable to all other algorithms considered for the problem \cite{Dries, gronlund},
  and performs better than all player solutions of \cite{sorensen2016}. Furthermore, we show that the player solutions were optimized with a different algorithm before being compared to the results from the KASS algorithm.
  The authors of \cite{sorensen2016} have acknowledged both findings. Finally, we show that in contrast to the claims in \cite{sorensen2016}, the players did not explore two different strategies.
  In fact, all the players followed the same strategy.
\end{abstract}
\begin{figure}[h]
  \centering
  \includegraphics[width=0.7\linewidth]{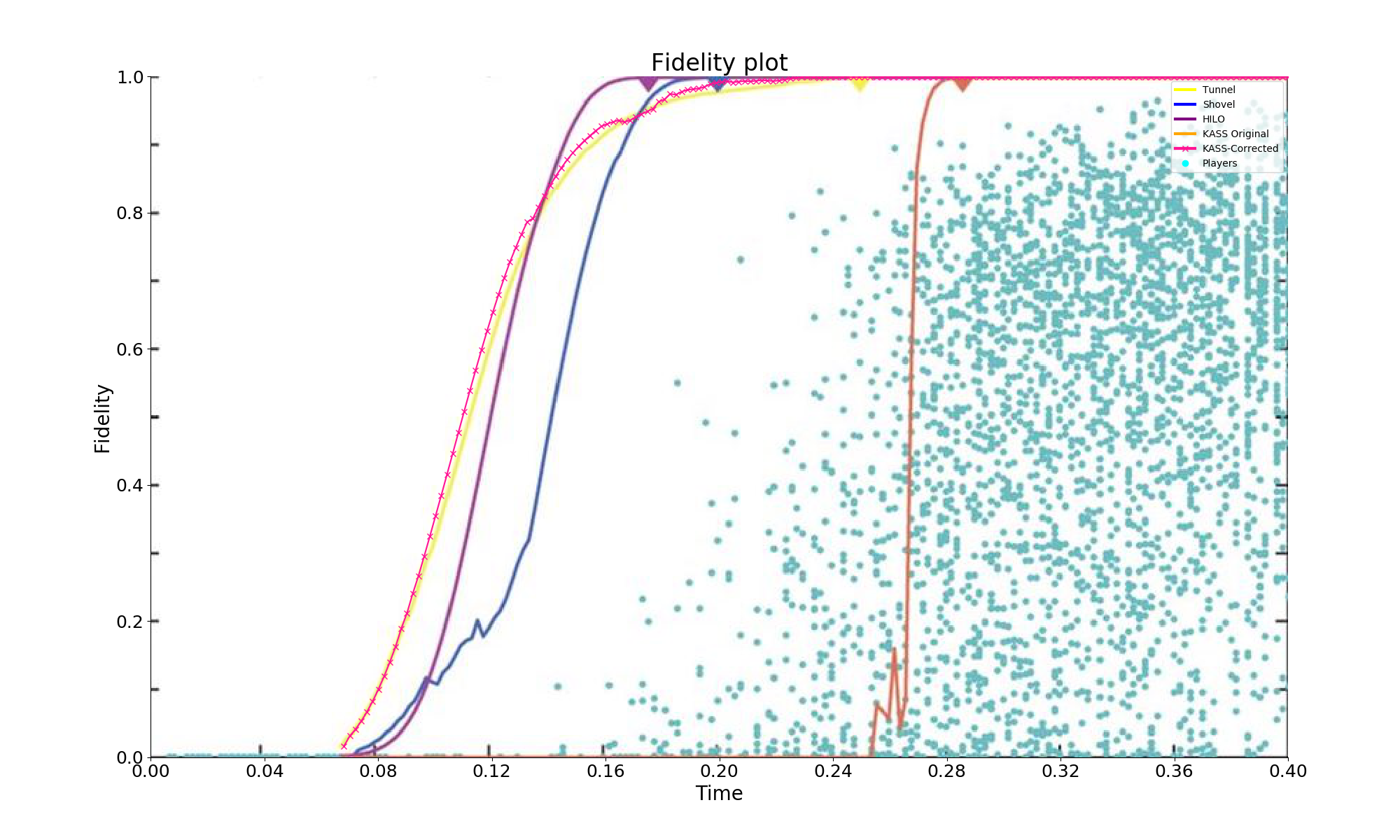}
  \caption{Fidelity for BringHomeWater algorithms from \cite{sorensen2016} (the colored curves) as a function of duration, along with player solutions (blue dots).
    Intuitively, the more to the left the better the solution. 
    Notice the massive difference between the KASS code from \cite{sorensen2016} to the right (orange curve),
    and the results from the exact same code (pink curve) with the corrected sign error (that was made on the basis of only 50 of the random starting positions and amplitudes from \cite{sorensen2016}).
    Note also, that the corrected KASS version (pink) is working extremely well despite that it has been optimized for
    code with the sign error and not the correct implementation, and that the algorithm outperforms all player solutions by a massive margin.
  }  
  \label{fig:all_compare}
\end{figure}
\clearpage

\section{Introduction and summary of results}
In the paper \emph{Exploring the quantum speed limit with computer games} \cite{sorensen2016} S{\o}rensen et al.  present and discuss results for a problem in quantum computing 
named {\it BringHomeWater}. They consider how well humans can solve the problem by playing a game where good game solutions map to good solutions for the BringHomeWater problem, and they compare the solutions found by human players to the solutions found by algorithms.
More specifically, they compare to solutions found by the so-called KASS algorithm, which is based on a standard algorithm from quantum optimization known as the  Krotov algorithm \cite{krotov}.
The claim of the paper is that players outperform KASS and that player intuition is essential for solving the problem effectively. The latter was emphasized by the authors by deriving a better algorithm than KASS from player solutions.
However, these claims were quickly disputed by D. Sels \cite{DriesarXiv, Dries}, who showed that another standard and very basic algorithm works much better than any player solution and algorithms derived from the player solutions.
D. Sels findings were later confirmed in \cite{gronlund} that also tested several other basic algorithms and showed that all these algorithms, including an implementation of the Krotov algorithm (also the basis of KASS),
give similar results and easily outperform all player solutions. Furthermore, any intuition found by players where easily found by just looking at the output from any one of the algorithms tested.
The code for these experiments was made publicly available at GitHub at the time of publication, and in attempt to to understand this large difference we asked the authors of \cite{sorensen2016} to explain the exact details of their implementation and provide their code, so we could compare and find the reason for the large difference in results. Instead of receiving any code or explanation of the approach in \cite{sorensen2016}, our results were subsequently disputed \cite{molmer} by a subset of the authors of \cite{sorensen2016}.
In particular, the superior performance of the Krotov algorithm of \cite{gronlund} was explained as being the result of employing an ``advanced arsenal of algorithms'' (like the ADAM method) to help the algorithm without any actual proof of this being the case.
Hence, they claimed that the results in \cite{Dries, DriesarXiv, gronlund} were somehow not in contradiction with the findings in \cite{sorensen2016}. However, the algorithm used by D. Sels \cite{Dries, DriesarXiv}
is a centuries old textbook method and the Krotov algorithm as implemented in \cite{gronlund} uses the most basic approach possible and relies only on standard off-the-shelf tools, that were available before the publication of \cite{sorensen2016}.

Finally, almost half a year after we asked, we have now received the actual (MatLab) code used in \cite{sorensen2016} from the authors, and we have investigated why the Krotov algorithm as considered in \cite{gronlund}
works so much better than the  Krotov based KASS algorithm used in \cite{sorensen2016}. In this note we report on the results of this investigation. 
Our main conclusions are as follows.
\begin{itemize}
\item We find that the Matlab code directly reproduce the KASS results of \cite{sorensen2016} and confirms that the KASS algorithm implementation of \cite{sorensen2016} does indeed  perform as poorly as shown.
  See Section \ref{sec:recreate}.  However, we have discovered that the inefficiency of the KASS algorithm is entirely caused by a simple sign flip error in the code, which we explain in Section \ref{sec:error}. 
  The error is directly visible from any analysis of the output of the KASS code provided. In particular, it is clear from simply plotting the output of the KASS algorithm.
  Note that no plots or basic visualizations of the output from the KASS algorithm or other similar analyses of the KASS algorithms failure are included in \cite{sorensen2016}.
\item If we fix the error (a simple erroneous sign in a hard-coded formula for a derivative) we obtain essentially
  the same result as all the other algorithms tested in \cite{gronlund}.
  In particular, the corrected algorithm outperforms all player solutions (Figure \ref{fig:all_compare}). This is explained in Section \ref{sec:fixed_alg}.
  Note that this means that the success of the Krotov algorithm of \cite{gronlund} is not caused by any ``advanced arsenal of algorithms'' (and in particularly not the ADAM method).
  It also intuitively means that the players do not contribute anything to the solution of the problem that basic implementations of standard algorithms cannot easily find autonomously.
  This is true now and it was also true at the time of the study of \cite{sorensen2016}.
\item We have discovered that player solutions cannot have been optimized with the same KASS algorithm as the random starting positions and amplitudes
  in direct contrast to what was stated in \cite{sorensen2016}. This is explained in Section \ref{sec:player}.
\item In \cite{sorensen2016} the players are credited for discovering and exploring two separate solution strategies with different properties.
  We show that this is not the case, irrespective of errors in the KASS algorithm. The players only consider and explore one strategy,
  and this is obvious from just plotting the player solutions (a plot not included in \cite{sorensen2016}). The credit for discovering the second solution strategy belongs solely to the Krotov algorithm.
  This is discussed in Section \ref{sec:yellow_blue}. It is also important to note that there is nothing special in finding these two strategies, as all algorithms tested easily find them both (\cite{gronlund, Dries, DriesarXiv}),
  including the KASS algorithm after we fixed the sign error. 
\item In \cite{sorensen2016} the final gained result is an improved algorithm, called \emph{HILO}, derived from the insight drawn from the two player strategies,
  and is credited to human problem-solving skills. However, as already mentioned the players are only responsible for finding one strategy, and all algorithms tested find both strategies.
  Furthermore, all the standard algorithms considered in \cite{Dries, DriesarXiv, gronlund} outperform the HILO algorithm.
  Hence, in contrast to the claims in \cite{sorensen2016}, it is clear that player does not provide any insights into the problem that standard algorithms does not.   
\end{itemize}
In summary, all the claims, observations and conclusions in the paper by Sorensen et al. \cite{sorensen2016} is based on the KASS algorithm implementation, which has a vital error that results in extremely poor results. 
Fixing the errors gives drastically different results, and verifies the results and conclusions from \cite{Dries, DriesarXiv, gronlund} while disproving all the claims and observations in \cite{sorensen2016} and the subsequent extra arguments given in \cite{molmer}.
Furthermore, the player solutions were optimized with a different algorithm before being compared with the results from the (erroneous) KASS algorithm (and not the same algorithm as claimed in \cite{sorensen2016}),
and the players are wrongly credited for discovering and exploring two distinct solution strategies when they clearly only consider one.
As a result of the first version of this note being shown to the authors of \cite{sorensen2016}, the dispute note \cite{molmer} has been withdrawn without giving any explanation of the issues we present.

\section{Recreating results for the KASS algorithm}
\label{sec:recreate}
\begin{figure}[ht]
  \begin{subfigure}[t]{0.48\textwidth}
    \centering
    \includegraphics[width=0.95\linewidth]{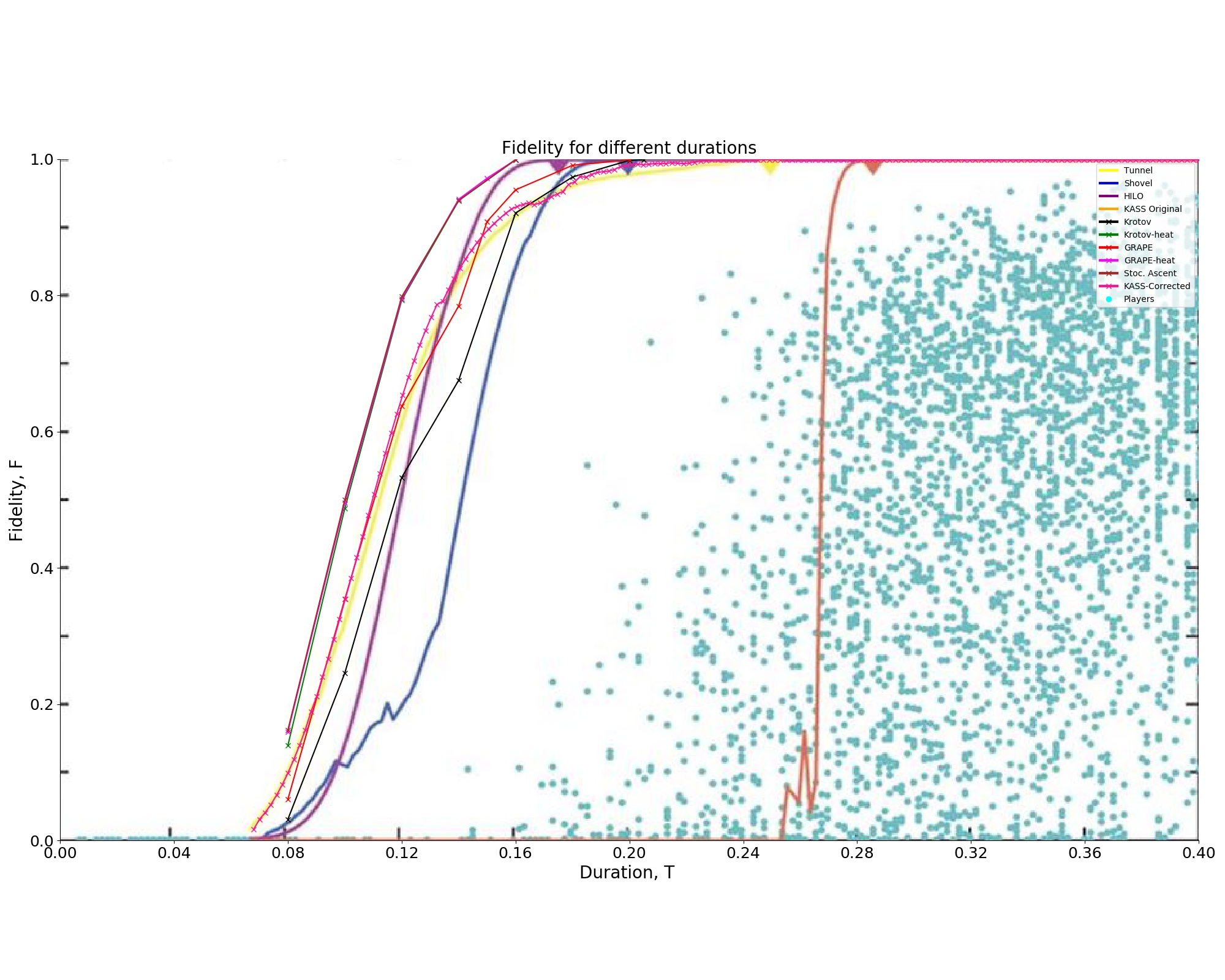}
    \caption{Comparison of all algorithms tested for BringHomeWater in \cite{sorensen2016, DriesarXiv, Dries, gronlund}, including the corrected KASS algorithm.}      
    \label{fig:all_overlay}
  \end{subfigure}%
  \hfill
  \begin{subfigure}[t]{0.48\textwidth}
    \centering
    \includegraphics[width = 0.9\linewidth]{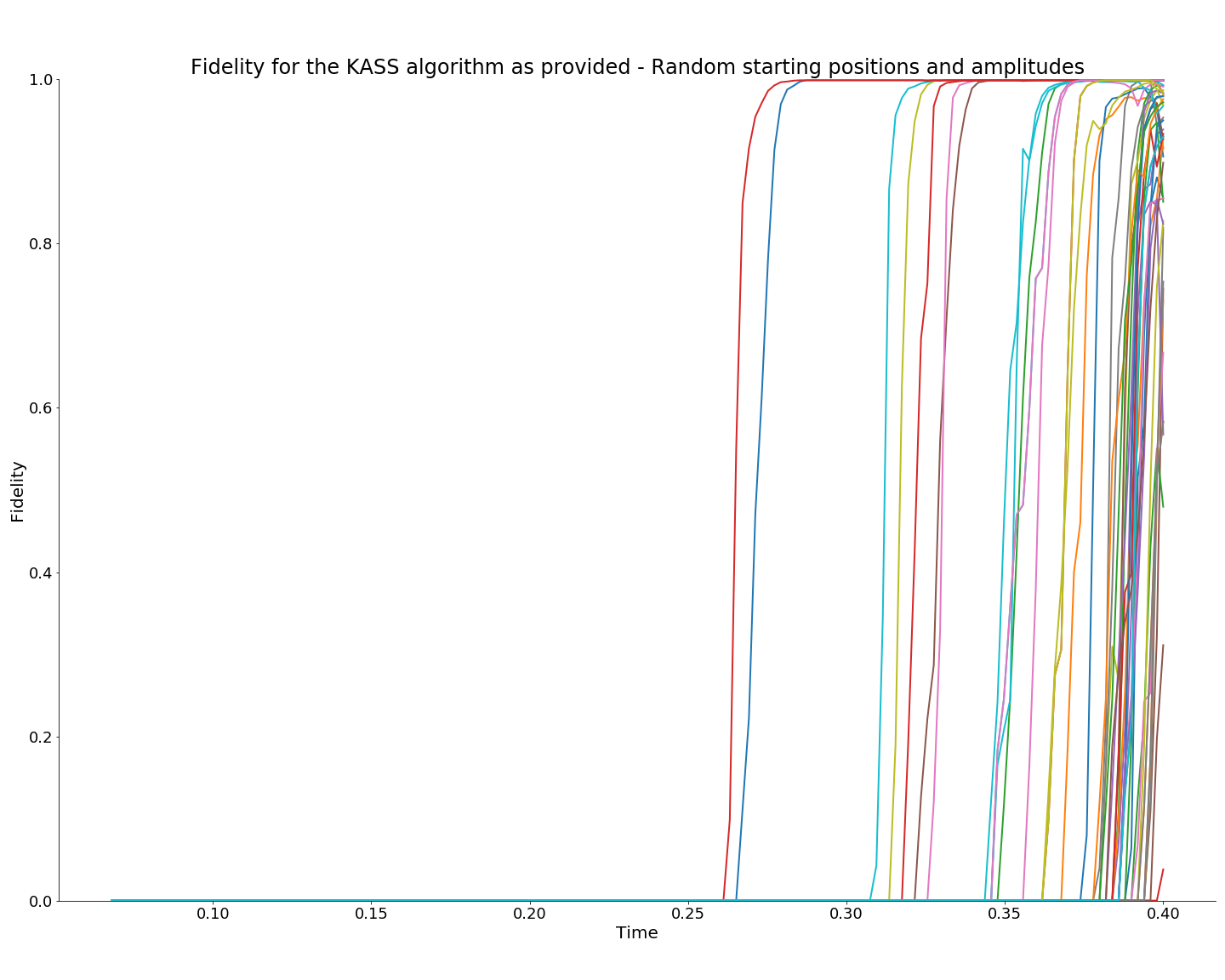}  
    \caption{Fidelity achieved by running the KASS code provided by the authors of \cite{sorensen2016} for a few of the random starting positions and amplitudes.
      The algorithm stops when the fidelity drops below 0.2.}  
    \label{fig:kass_as_given}
  \end{subfigure}
  \caption{Results for the KASS algorithm}
\end{figure}

By using the exact code used in \cite{sorensen2016} provided by the authors, we have reproduced the results for the KASS algorithm using a small subset of the random starting positions and amplitudes from \cite{sorensen2016}, and the results are shown in Figure \ref{fig:kass_as_given}.
As we can see, the plot is in perfect agreement with the result presented for the KASS algorithm in \cite{sorensen2016} (Figure \ref{fig:all_compare}).

We have shown the results for different algorithm implementations in Figure \ref{fig:all_overlay} taken from \cite{gronlund}, where we have included the corrected KASS implementation for comparison.
Here it is obvious, as it was shown in \cite{gronlund}, that the only algorithm that does not work for the problem is the KASS implementation from \cite{sorensen2016}.
The main differences between the Krotov based KASS algorithm in \cite{sorensen2016} and the implementation of the Krotov algorithm in \cite{gronlund} are as follows
\begin{itemize}
\item \textbf{Small change of problem parameters}. The parameters for the BringHomeWater problem are not given in \cite{sorensen2016}.  In \cite{gronlund} we used the parameters specified in \cite{Dries} that tried to infer them from the actual game itself.
  In the note provided by the authors of \cite{sorensen2016} along with their code, we are given the actual parameters and there are some small differences. None of these differences have any consequence for the overall conclusion of our investigations in \cite{gronlund}.
\item \textbf{Initial starting positions and amplitudes:} In \cite{gronlund}, we tested uniform random positions and amplitudes (the most naive method possible) and an
  approximation of the starting positions and amplitudes used in \cite{sorensen2016} (as they were not fully specified).
  The starting positions and amplitudes used in \cite{sorensen2016} is given by an involved linear combination of sine functions,  and we will not discuss their inner workings here.
\item \textbf{Sweeps:}  In \cite{gronlund}, the Krotov algorithm is tested on each duration separately i.e.  the algorithm is run independently for each fixed duration.
  In \cite{sorensen2016}, the authors introduce a \emph{sweep} that works as follows. The sweep starts by applying the Krotov algorithm at high duration (0.4), where the problem is very easy, on a random set of starting positions and amplitudes.
  The Krotov algorithm is then run for 800 steps or until a fidelity score of 0.999 is achieved.
  The last values for the positions and amplitudes found by the Krotov algorithm are then used as input for the same Krotov algorithm on a slightly smaller duration.
  Repeating this process until the duration becomes very low (0.068) is known as a sweep, and gives a score for each duration considered in the sweep. 
  \item \textbf{Learning rates}. In \cite{gronlund} we use an off-the-shelf learning rate controller, known as ADAM \cite{kingma2014adam},  that was available at the time of the study, where \cite{sorensen2016} chose to implement their own heuristic for controlling the learning rates to improve their KASS algorithm.
\end{itemize}
To see if these differences may explain the large discrepancy in the results between \cite{gronlund} and \cite{sorensen2016},
we updated the implementation used in \cite{gronlund} with the provided parameters, and replaced the ADAM learning rate controller with the naivest method possible of using fixed learning rates. We have implemented the sweeps and tested the same starting positions and amplitudes as considered in \cite{sorensen2016}.
The result of these experiments matches the results for the Krotov algorithm as given in \cite{gronlund}.
Hence, none of these four differences are the reason for the discrepancy, in particular it has nothing to do with ADAM or any other ``advanced arsenal of algorithms'' as claimed as the reason (without any evidence) in \cite{molmer}.

\section{Finding an error}
\label{sec:error}
\begin{figure}[ht]
  \begin{subfigure}[t]{0.45\textwidth}
    \centering
    \includegraphics[width=0.85\linewidth]{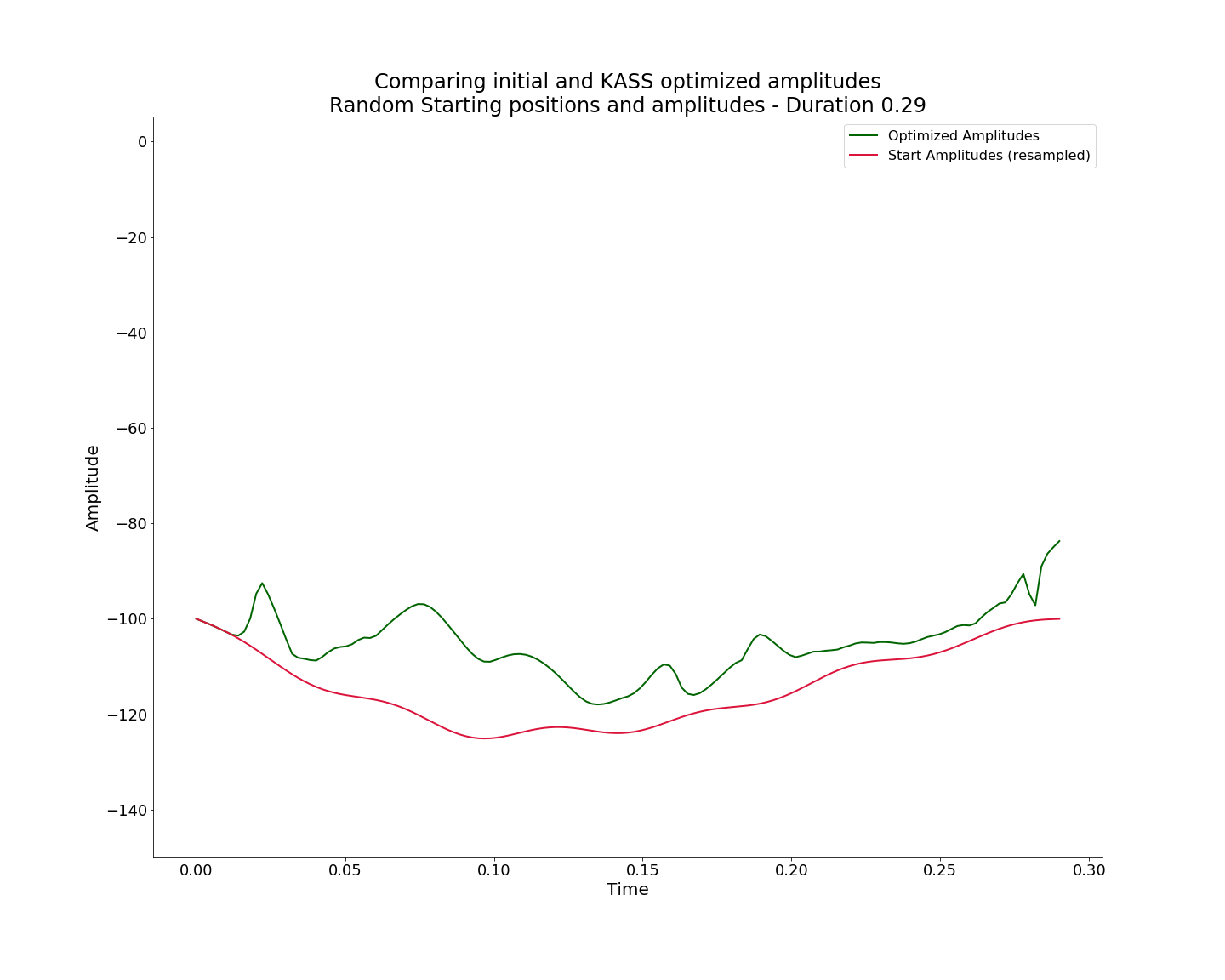}  
    \caption{The starting amplitude and the amplitude as optimized by the KASS algorithm for a single starting position and amplitude at duration 0.29,
      the last duration with high fidelity. The red line is the initial amplitude and the green is the optimized amplitudes. }      
    \label{fig:amp_decay1}
  \end{subfigure}%
  \hfill
  \begin{subfigure}[t]{0.45\textwidth}
    \centering
    \includegraphics[width=0.85\linewidth]{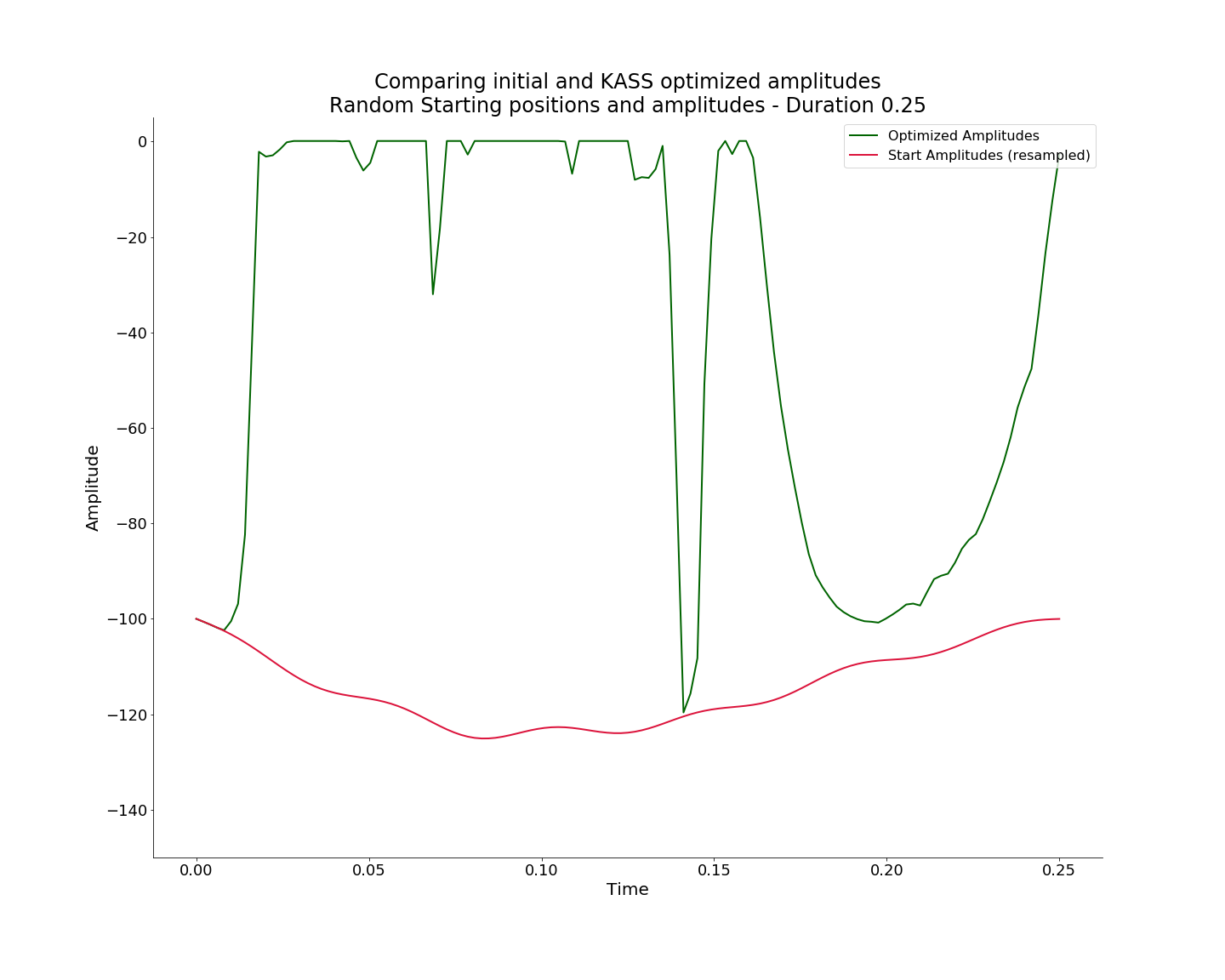}  
    \caption{The starting amplitude and the amplitude as optimized by the KASS algorithm for a single starting position and amplitude at duration 0.25.
      The red line is the initial amplitude and the green is the optimized amplitudes.}      
    \label{fig:amp_decay2}
  \end{subfigure}\\
  \begin{subfigure}[t]{0.45\textwidth}
    \centering
    \includegraphics[width=0.85\linewidth]{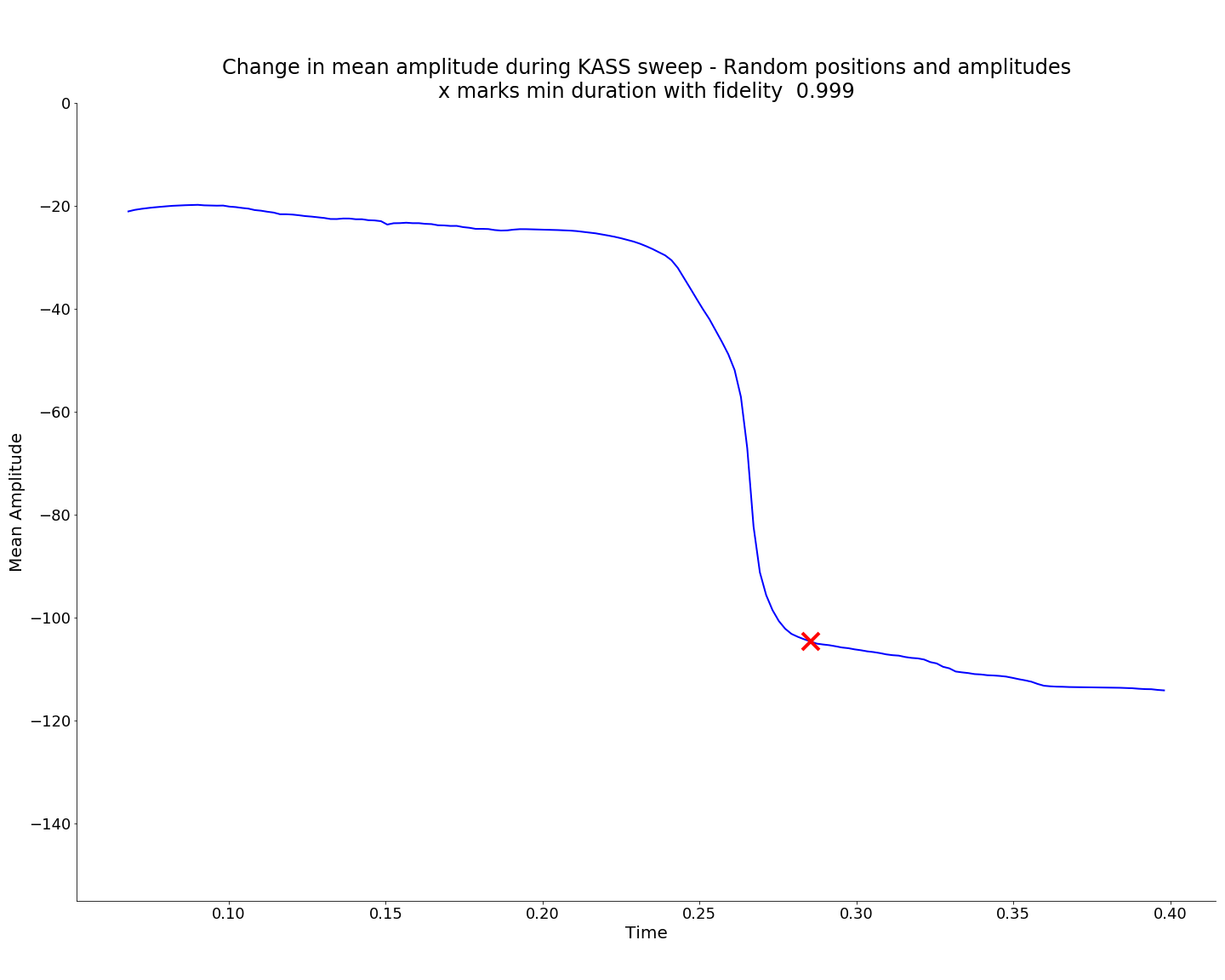}   
    \caption{The mean amplitude as a function of the duration. The red cross marks the last duration where 0.999 in fidelity was achieved.
      Note that the plot is read from right to left in time, as the sweep starts from duration 0.4 and decrease duration in each step.}      
    \label{fig:amplitude}
  \end{subfigure}%
  \hfill
  \begin{subfigure}[t]{0.45\textwidth}
    \centering
    \includegraphics[width=0.85\linewidth]{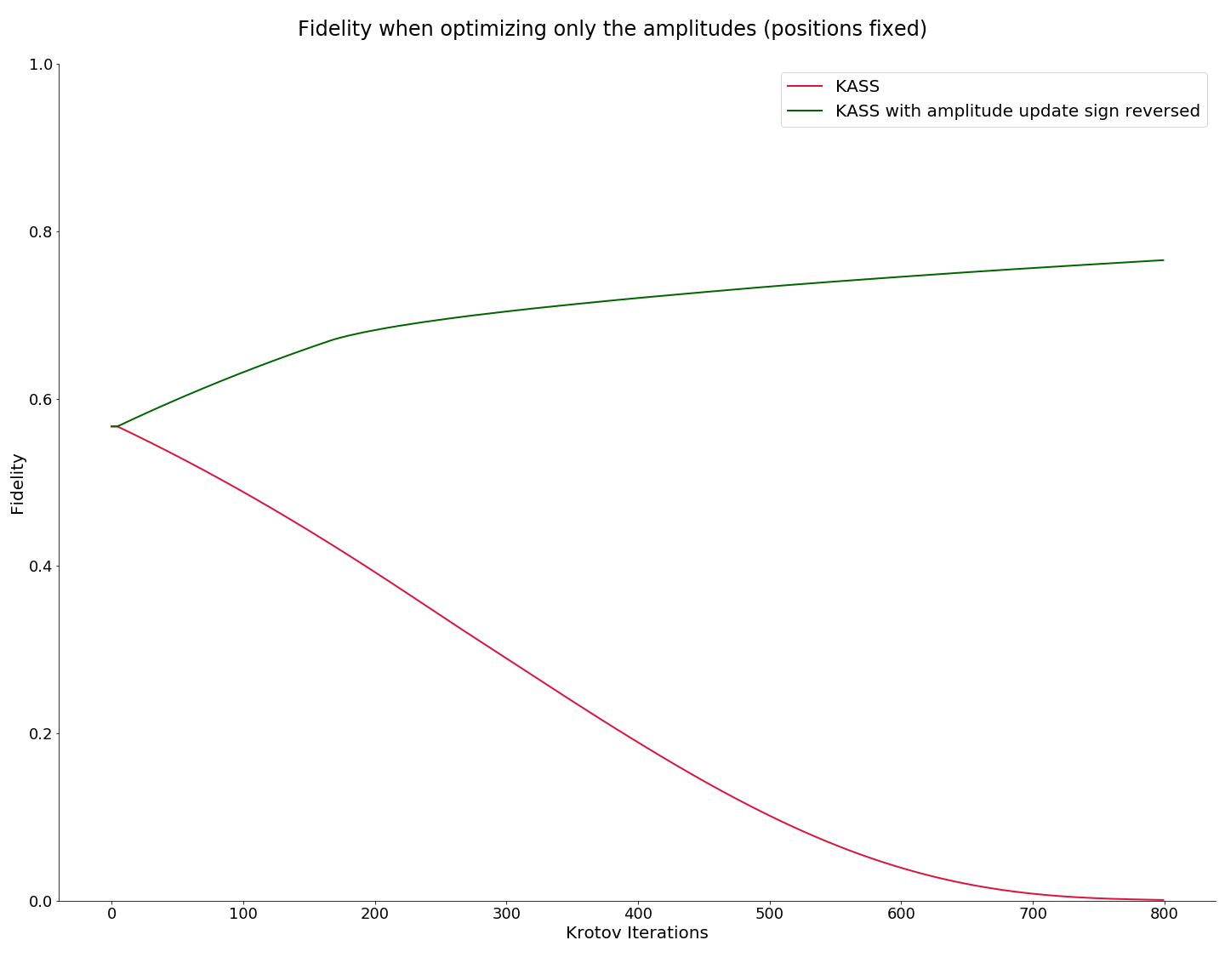}  
    \caption{Fidelity as a function of the number of Krotov iterations for the KASS implementation in \cite{sorensen2016}.
      The red line is the fidelity achieved with the original code, the green line is the fidelity achieved with the sign for the amplitude updates reversed.}
    \label{fig:fixed_pos}
  \end{subfigure}%
  \caption{Visualizing the failure of the KASS Algorithm}  
\end{figure}

We looked over the Matlab code provided by the authors of \cite{sorensen2016}, computed a few basic statistics when running the KASS algorithm, and visualized the results.
The thing that stands out when using the KASS algorithm as provided by \cite{sorensen2016}, is that the amplitudes keep getting closer and closer to zero as more and more iterations of the algorithm is performed.
In Figure \ref{fig:amp_decay1}, and \ref{fig:amp_decay2} we have showed the starting amplitudes, and the amplitudes that result from running the KASS code provided.
From this simple plot of the output of the of KASS algorithm, it is clear what the algorithm is doing is wrong. 
In laymans terms, the algorithm turns the tweezer controlled by the algorithm more and more off in each iteration.  This makes it impossible to move the atom and the result is low fidelity solutions.
This is in stark contrast to what we observe using the code used in \cite{gronlund}, where all experiments show that the Krotov algorithm essentially always turn the tweezer amplitude to maximum power (-150)\footnote{-150 is the maximum power output and 0 is the minimal}.
The same happens in the GRAPE \cite{grape} algorithm, and the Stochastic Ascent \cite{Dries} algorithm tested for the problem, and setting the amplitude to the maximum output power also fits perfectly with theory \cite{Dries}.

To investigate the situation more, we have showed the development of the mean amplitude for a single set of starting positions and amplitudes in Figure \ref{fig:amplitude}.
Here we observe that the amplitudes are initially slowly dropping towards zero, but when the algorithm can no longer find a  solution with fidelity 0.999, it starts to move the amplitude rapidly towards zero, explaining the rapid descent from a high fidelity solution to essentially zero fidelity with only a small reduction
in the duration time considered. We note that the results from the KASS algorithm \cite{sorensen2016} are the only results that shows this behavior as all other algorithms slowly decrease in fidelity as the duration given is decreased, including the results from starting with from player solutions
that was claimed to be optimized with the same algorithm in \cite{sorensen2016} (see Figure \ref{fig:all_overlay}).
The reason the KASS algorithm actually achieves a fidelity of 0.999 until a duration of 0.29, is that the algorithm has used relatively few updates to positions and amplitudes to get fidelity 0.999 for the high durations, and improving the positions help more than decreasing the amplitude initially.
This may have been further boosted by the fact that the KASS code as provided by the authors of \cite{sorensen2016}, caps the learning rates which ensure that not too large updates are made to the amplitudes.

To test our findings, we made an additional experiment where we change the KASS implementation such that it only updates the amplitudes. We initialize the algorithm with very good starting positions and reasonably large amplitudes. We then test the algorithm with a version for updating the amplitude as it was provided, and a version where we flip the sign of the amplitude update.
The results are shown in Figure \ref{fig:fixed_pos}. Here it is clear that if the algorithm, as provided by the authors of \cite{sorensen2016},
is only allowed to optimize amplitudes the fidelity quickly drops. However, if we change the sign of the amplitude update, the fidelity  improves as it should in a local search algorithm as Krotov. 

Searching for the reason for this behavior, we found a sign error in the formula for computing the derivative in regard to the amplitudes.
The technical derivation is given in the Appendix. It is not a numerical error, but an error in a hard-coded formula for the derivative, and it means that the KASS implementation from \cite{sorensen2016} tries to find the worst possible solution for the amplitudes instead of the best. We note that any basic analysis of the output of the KASS algorithm would reveal this obvious flaw in the implementation.

\section{Results with a correct KASS implementation}
\label{sec:fixed_alg}

\begin{figure}[ht]
  \begin{subfigure}[t]{0.48\textwidth}
    \centering
    \includegraphics[width = 0.9\linewidth]{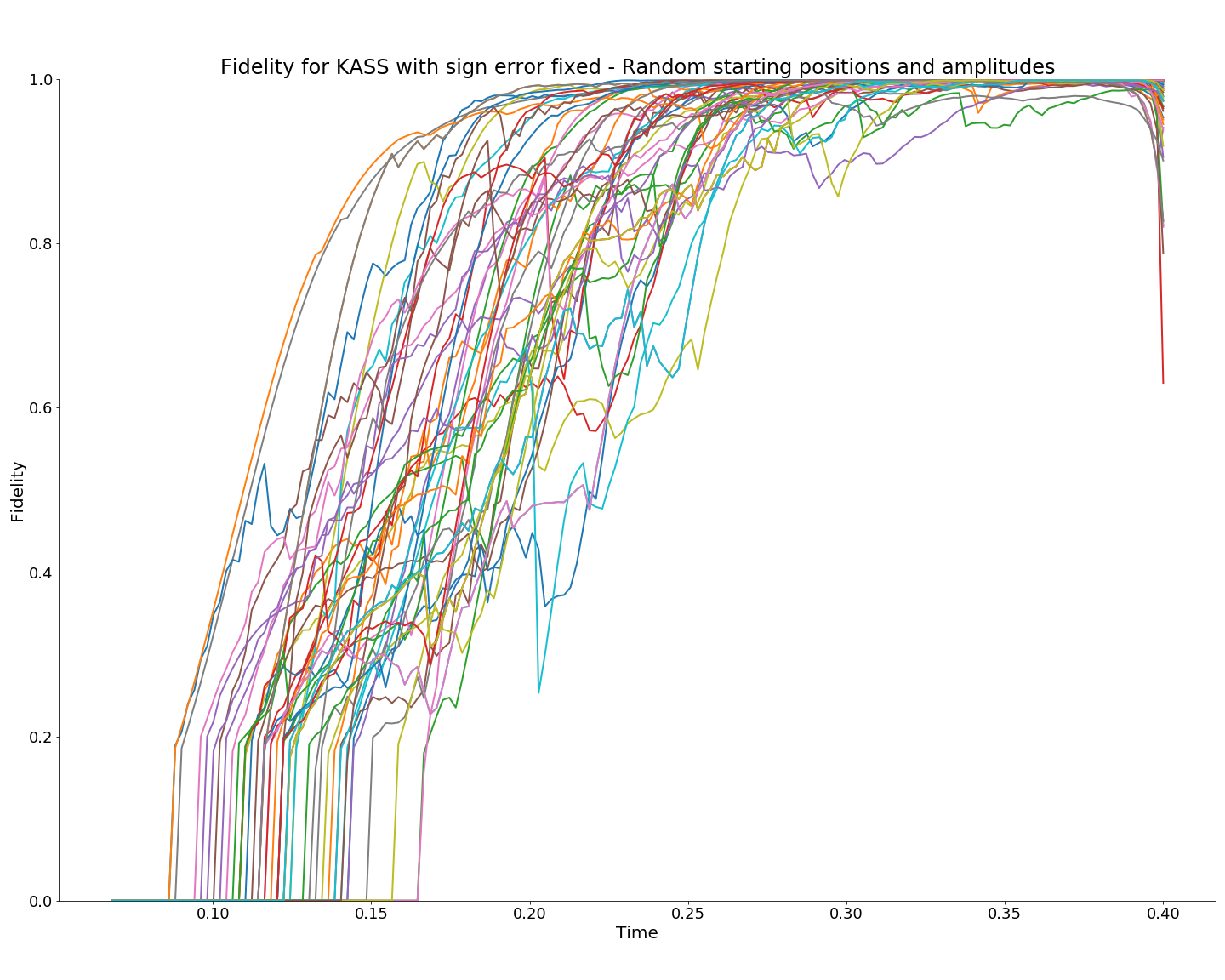}  
    \caption{Fidelity achieved by the correct KASS algorithm (correct sign for the derivative of the amplitude), for a few of the different starting positions and amplitudes from \cite{sorensen2016}. 
      As before the algorithm stops when then the fidelity drops below 0.2, explaining the drastic drop that finish off each result.}
    \label{fig:kass_fixed}
  \end{subfigure}%
  \hfill
  \begin{subfigure}[t]{0.48\textwidth}
    \centering
    \includegraphics[width=0.9\linewidth]{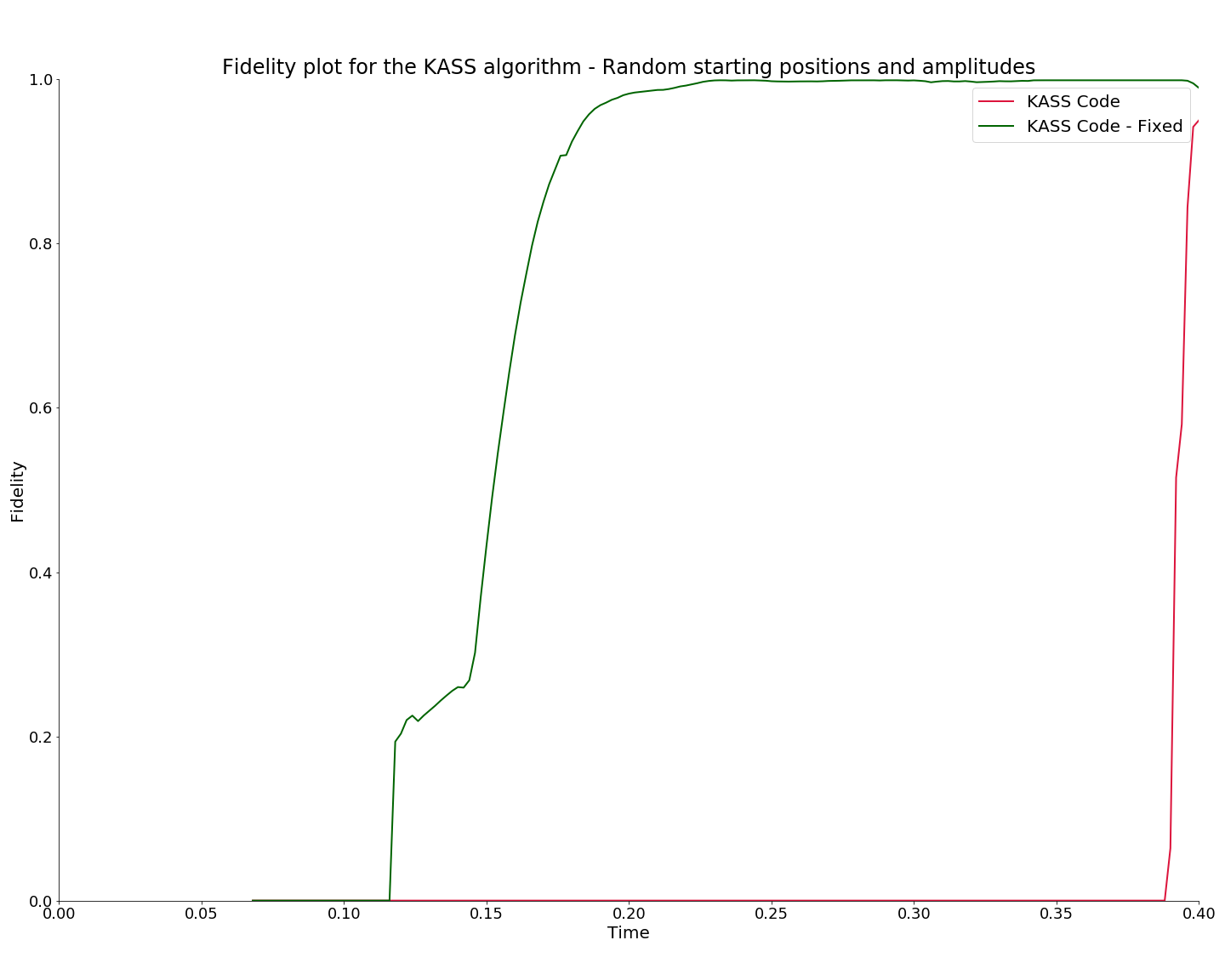}  
    \caption{Comparing the fidelity achieved for the KASS algorithm with and without the sign error on a single pair of random starting positions and amplitudes.
      The red line shows the poor performance of the original KASS algorithm, and the green line shows the greatly improved performance of the KASS algorithm with the error fixed.}
    \label{fig:compare_994}
  \end{subfigure}%
  \caption{Results for a correct implementation of the KASS algorithm.}
  \label{fig:fix_error}
\end{figure}

We remedied this sign error, and ran the otherwise exact same KASS algorithm implementation on the same starting positions and amplitudes as used for Figure \ref{fig:kass_as_given}, and the results are shown in Figure \ref{fig:kass_fixed}.
As this figure shows, using only a few of the random starting positions and amplitudes considered in \cite{sorensen2016}, we get a result that looks just like the result given in \cite{gronlund},
and this is with an algorithm that is not optimized with hyper-parameters and learning rates for an actually correct implementation.
To underline the point, we have shown the difference between running the KASS algorithm as provided, with and without the sign error, on a single starting position and amplitude pair in Figure \ref{fig:compare_994},
and we can see there is a world of difference between the results.
The sign error also explains why the fidelity achieved using the erroneous KASS algorithm drops from a high value to essentially zero in a very short duration span, as it explains why the algorithm simply turns the tweezer off, making any transport of the atom impossible.

In \cite{sorensen2016}, the authors state that they have tested 2400 different random starting positions and amplitudes (their code and note reveals that this number is actually 1961) each taking 6 hours of computation time, and giving a total of $7.4 \times 10^8$ iterations of the Krotov algorithm, and use this as profound evidence
of more general nature that computers, and particular the Krotov algorithm, fails completely for the BringHomeWater problem. A problem that human players and their intuition can help to solve.
The truth is that all this computation time, and the very many iterations has been used to turn off the tweezer that is supposed to transport the atom. The paper also states that they \emph{''tried different seeding and optimization strategies and the most successful was a Krotov algorithm using sinusoidal seed functions and sweeps over the total duration (KASS)''}. Now, changing seeding strategy (random starting positions and amplitudes) may not help the KASS algorithm as implemented in \cite{sorensen2016} unless the amplitudes are increased. However, it is very clear that no optimization strategy, other than applying an erroneous algorithm, have been correctly tested by the authors of \cite{sorensen2016} as KASS is the only algorithm we have seen with such poor performance \cite{gronlund, Dries}. 

The results obtained with the KASS algorithm, and the arguments and conclusions based on it in \cite{sorensen2016} are therefore completely voided and in fact incorrect. The result confirms the conclusions made in \cite{gronlund, Dries}, and disprove the claims and dispute of these results from \cite{molmer} (that has now been withdrawn).
We note that using our own Krotov implementation from \cite{gronlund} initialized with the same starting positions and amplitudes as the KASS algorithm,
we easily get 0.999 fidelity at duration 0.18 without any tuning, and using this optimized solution also gives 0.999 in fidelity when used at duration 0.176 matching the result of the HILO algorithm out of the box. Performing the same basic analysis performed in \cite{gronlund} on the output of the algorithm leads to the same results as in \cite{gronlund}.
We get similar results when we use the player solutions to start the algorithm from instead.

\section{Testing Player and HILO Solutions.}
\label{sec:player}
\begin{figure}[ht]
  \begin{subfigure}[t]{0.48\textwidth}
    \centering
    \includegraphics[width=0.88\linewidth]{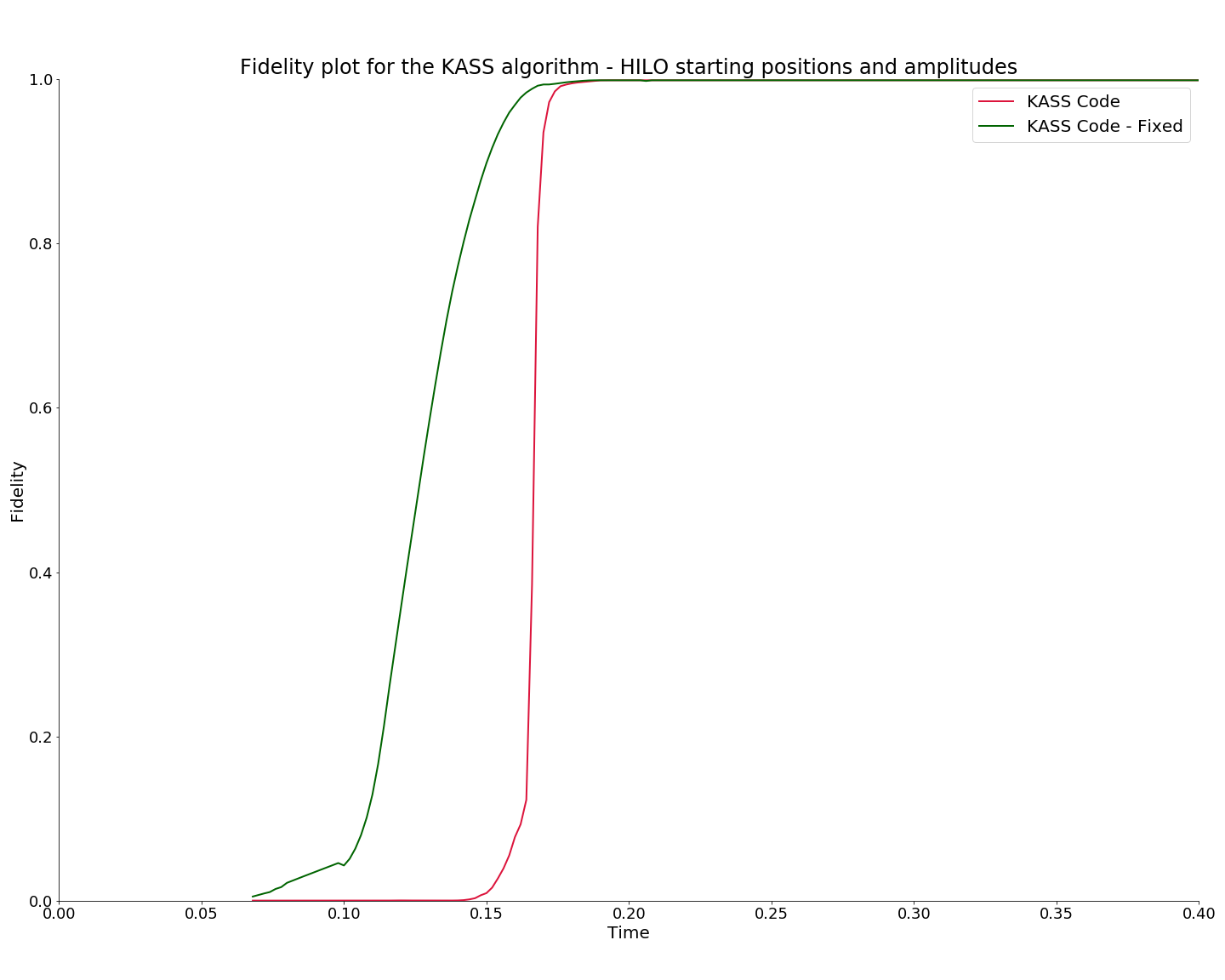}  
    \caption{Comparing KASS results with HILO starting positions and amplitudes. The red line shows the original algorithm and the green shows the corrected version.}
    \label{fig:hilo_compare}
  \end{subfigure}%
  \hfill
  \begin{subfigure}[t]{0.48\textwidth}
    \centering
    \includegraphics[width=0.88\linewidth]{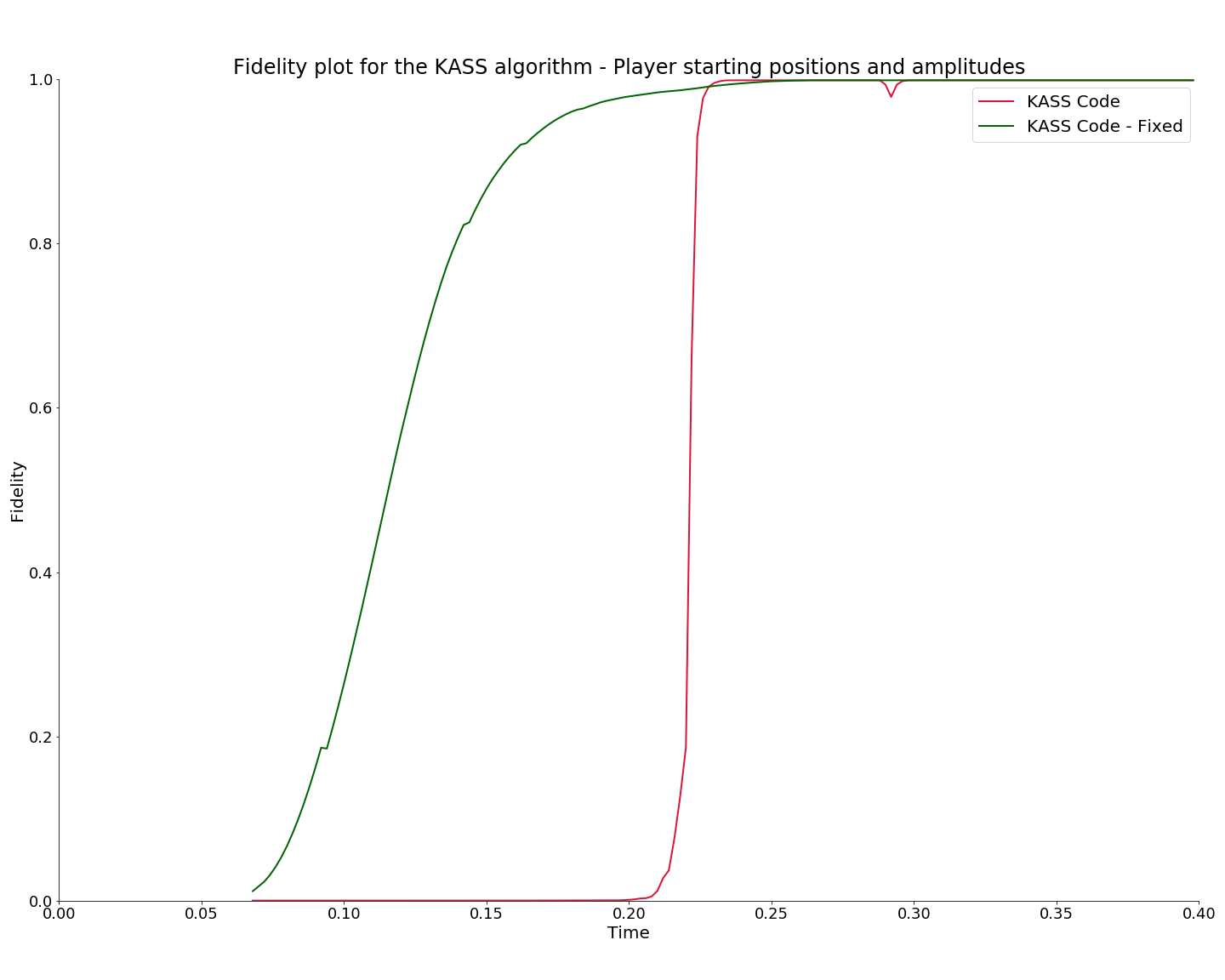}
    \caption{Comparing KASS results with player starting positions and amplitudes. The red line shows the original algorithm and the green shows the corrected version.}    
    \label{fig:yellow_compare}
  \end{subfigure}%
  \\
  \begin{subfigure}[t]{0.48\textwidth}
    \centering
    \includegraphics[width=0.88\linewidth]{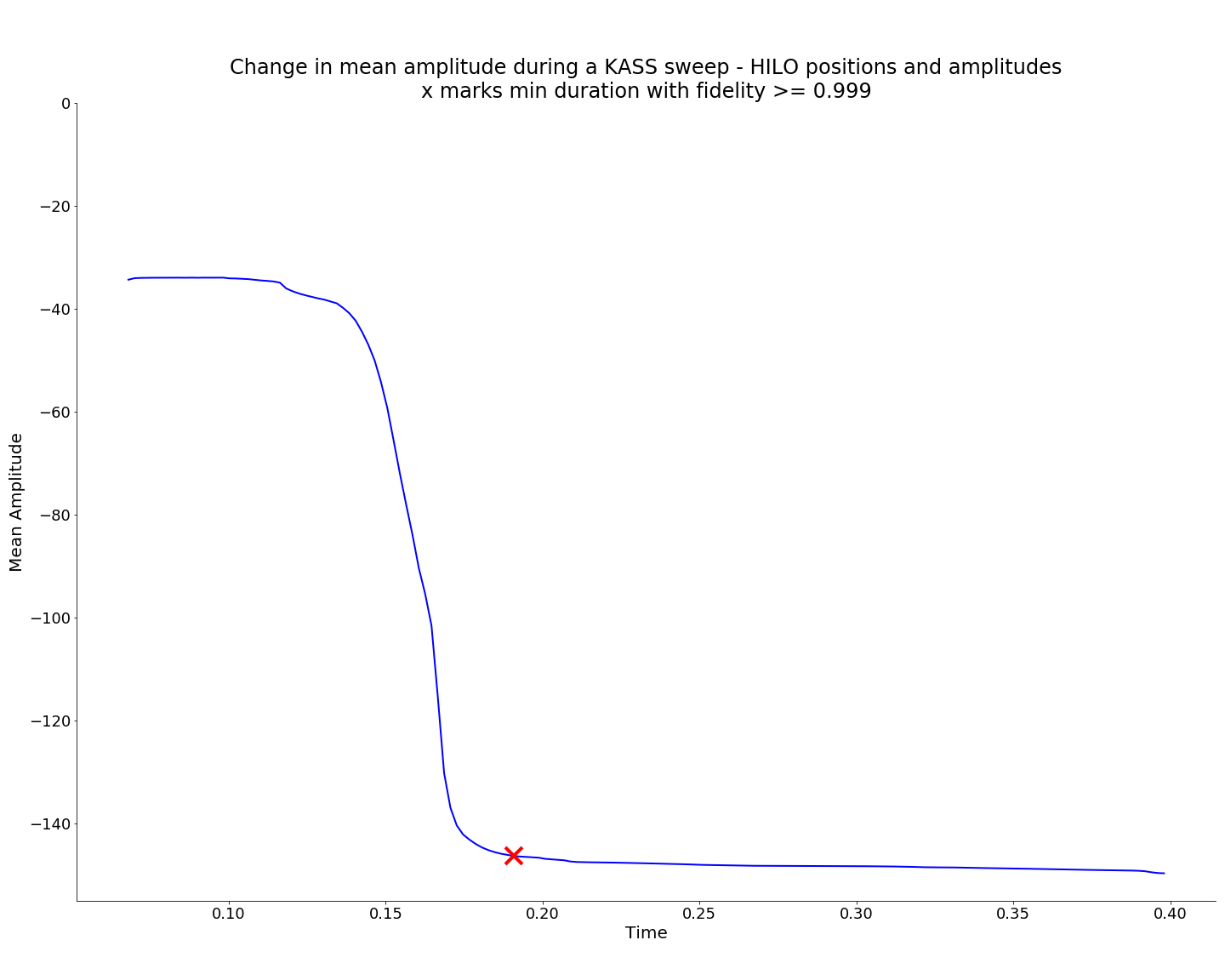}  
    \caption{The mean amplitude as a function of the duration for HILO starting positions and amplitudes.
      The red cross marks the last duration with 0.999 in fidelity.}
    \label{fig:hilo_amp_mean}
  \end{subfigure}%
  \hfill
  \begin{subfigure}[t]{0.48\textwidth}
    \centering
    \includegraphics[width=0.88\linewidth]{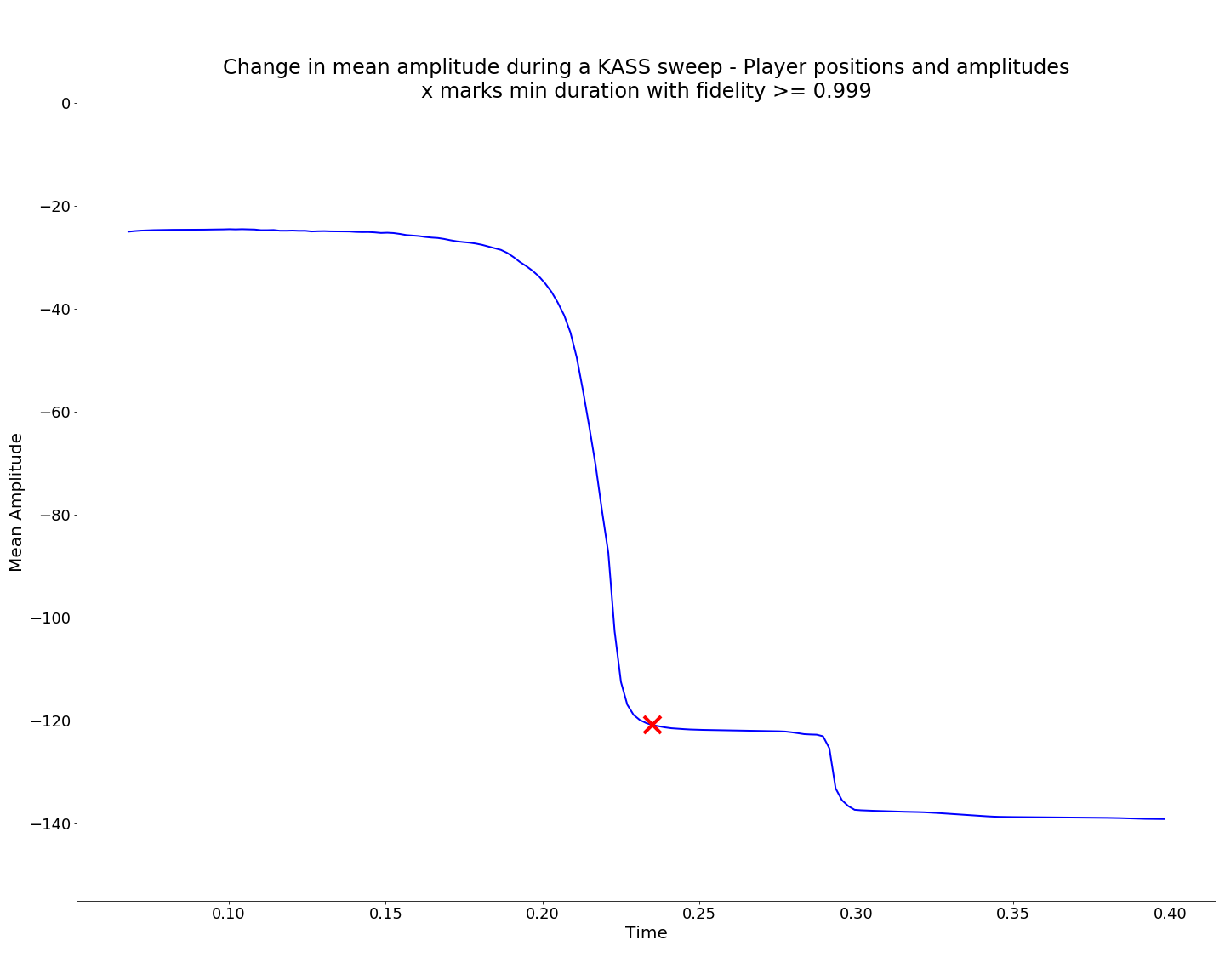}
    \caption{The mean amplitude as a function of the duration for the player solution.
      The red cross marks the last duration with 0.999 in fidelity.}    
    \label{fig:yellow_chop_mean}
  \end{subfigure}%
  \caption{Testing HILO and player starting values with original and error corrected KASS algorithm.}
  \label{fig:hilo_chop}
\end{figure}

Since the KASS algorithm as implemented in \cite{sorensen2016} has a sign error, you may wonder why the authors of \cite{sorensen2016} are able to achieve good results by starting the same KASS algorithm with the positions and amplitudes found by players,
and the positions and amplitudes found with the HILO algorithm developed in \cite{sorensen2016}. Secondly, why does the fidelity curves decay naturally when duration is decreased and the problem becomes harder when starting from these player solutions,
when this was clearly not the case for the KASS algorithm with random initialization (Figure 2 \cite{sorensen2016}).
Well, in the HILO algorithm the amplitude is fixed at the highest absolute value, -150, and if the algorithm does not update it then that works independent of the sign error in the KASS algorithm.
We have shown the results of running the KASS algorithm on the HILO starting positions both with and without the sign error, to see if we can reproduce the result from \cite{sorensen2016}.
This is shown in Figure \ref{fig:hilo_compare}, and in Figure \ref{fig:hilo_amp_mean} we show the development in the mean amplitude when using the KASS algorithm with the sign error.
Here we can see that the KASS algorithm (with the sign error) again slowly decrease the amplitude until the fidelity of 0.999 is no longer possible, and then the amplitude rapidly drops.
Exactly the same behavior as we saw earlier, a very quick drop from high fidelity to zero caused by the fact that the tweezer is turned off by the KASS algorithm.
If we run the KASS algorithm without the error, the HILO starting positions and amplitudes works better and the fidelity decays naturally, just like the behavior of all other algorithms tested on the problem (Figure \ref{fig:all_compare}), and the curve looks like the HILO curve in \cite{sorensen2016}.

But what about the player solutions.  When looking at the optimized positions and amplitudes found with player solutions (Extended Figure 3 of \cite{sorensen2016}), we notice that the amplitude for both optimized player starting positions and amplitudes shown, max out at -150
for large parts of the optimized solution at low duration. This is very hard to understand in light of the fact that there is a sign error in the amplitude derivative in the KASS code provided by the authors of \cite{sorensen2016} that should cause the amplitude to decrease,
as it has for all other starting positions and amplitudes we have tested (including the HILO starting positions and amplitudes from \cite{sorensen2016}).
We have shown the result of running the same experiment as we did for HILO above where we have initialized the algorithm with a player solution (from the yellow group) instead in \ref{fig:yellow_compare} and \ref{fig:yellow_chop_mean}.
Again, we see the same result of a fast drop. The intriguing double drop seems to fit with the small drop in the fidelity we see in the fidelity plot.
The only conclusion is that the player solutions have been optimized with a different algorithm, in contrast to what was stated in \cite{sorensen2016}. 

After our first version of this note were shown to the authors of \cite{sorensen2016}, they acknowledged this second critical error as well and told us what they really did.
In short, the player solutions were optimized with a different version of the KASS algorithm that had mitigated the error in the amplitude derivative, by making the learning rate for the amplitude so small that the 
algorithm could only makes tiny changes to the amplitudes, and therefore returns the same amplitude values as it was given as input. Another change, was a 25 percent increase in the number of iterations of the Krotov algorithm used in each step of the sweep. 
We conclude that the reason the player solutions get higher fidelity  than the random starting positions and amplitudes tested in \cite{sorensen2016},
is simply that the player solutions start with much larger amplitudes. Stated differently, the players learned  to turn the tweezer on, which was critical for success with the KASS implementation by \cite{sorensen2016}.
That, of course, has absolutely no relevance, as the optimization algorithm is incorrectly implemented to make the worst possible solution for the amplitude,
and the initial random amplitudes considered in \cite{sorensen2016} are so small that no good solutions for low durations are possible.

\section{Analysis of the player solutions}
\label{sec:yellow_blue}
\begin{figure}[ht]
  \begin{subfigure}[t]{0.48\textwidth}
    \centering
    \includegraphics[width=0.9\linewidth]{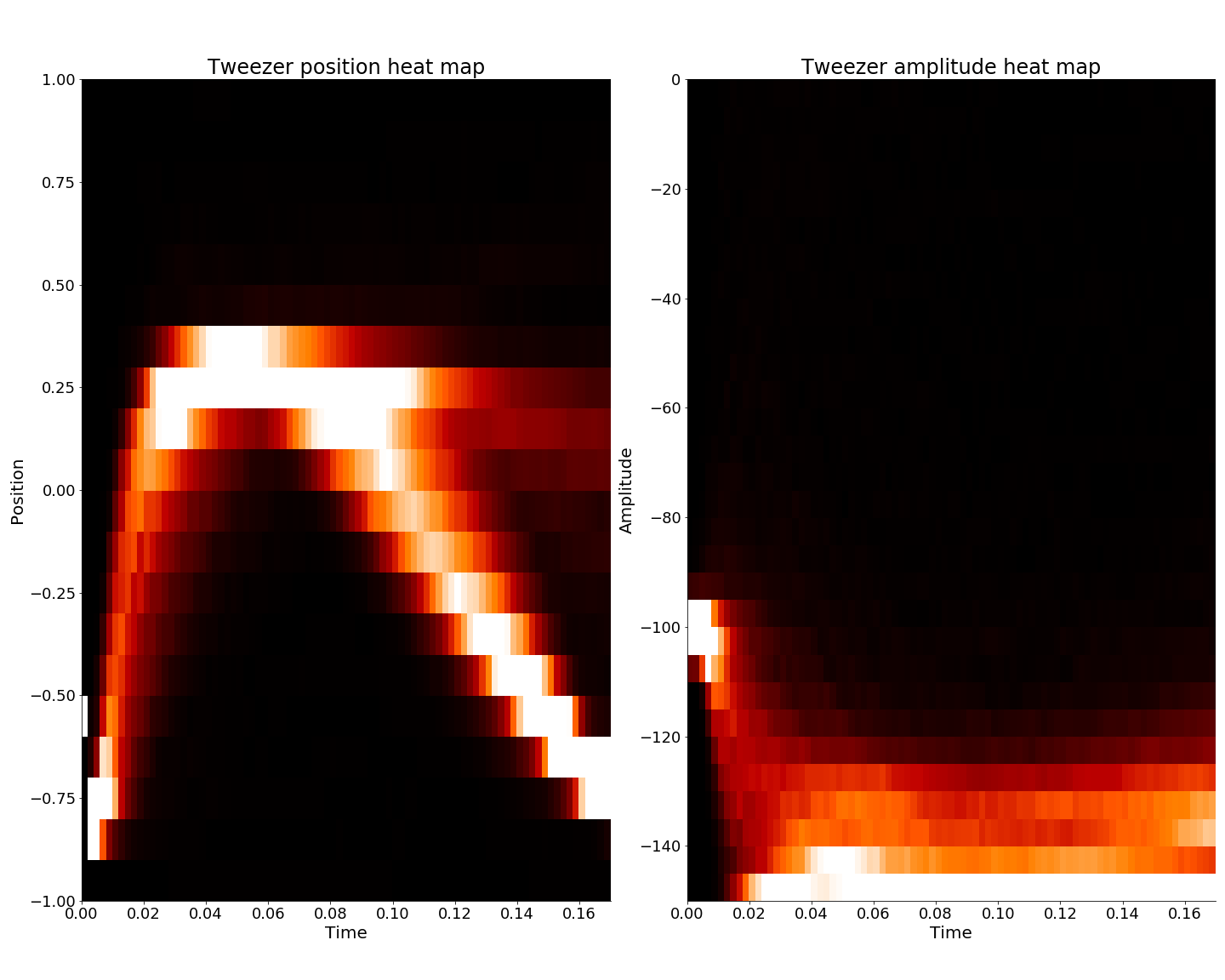}  
    \caption{Heat map of player solutions re-sampled to duration 0.17 (85 steps), but otherwise unchanged. 
      Each cell color displays the number of times a player solution position/amplitude was in the given interval at the given time step.
      The brighter the color the higher the count. }
    \label{fig:player_heat}
  \end{subfigure}%
  \hfill
  \begin{subfigure}[t]{0.48\textwidth}
    \centering
    \includegraphics[width=0.9\linewidth]{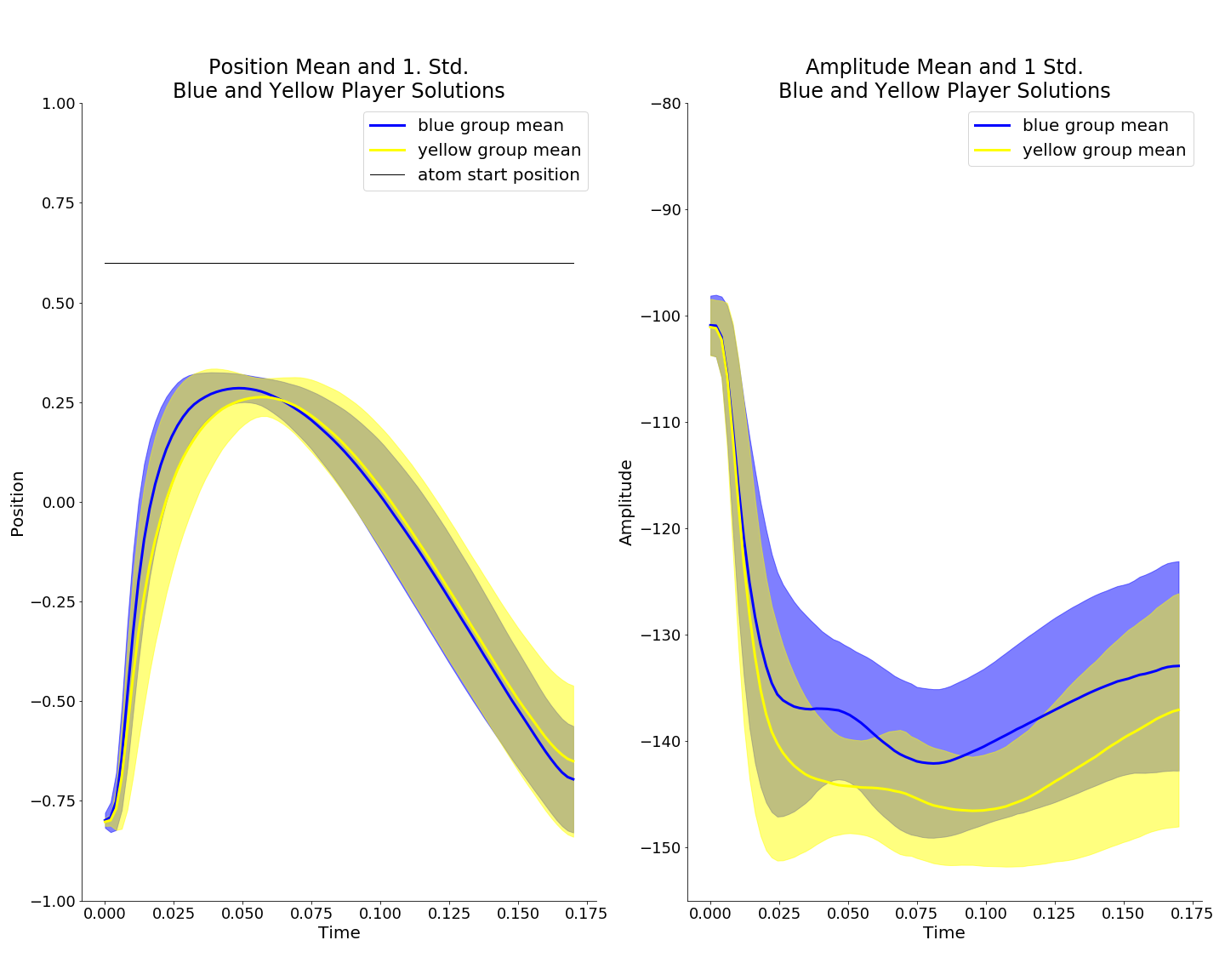}
    \caption{The mean and one standard deviation of the player solutions in the blue and the yellow clan respectively.
      The left plot shows the positions and the right plot shows the amplitudes.
      The black line in the left plot shows the initial position of the atom.}      
    \label{fig:yellow_blue_compare_mean_std}
  \end{subfigure}%
  \caption{Player solutions visualized}
\end{figure}
In \cite{sorensen2016} the authors discuss two different solutions strategies (see Figure 3 in \cite{sorensen2016}) they denote \emph{Tunneling} (yellow color in plots)
and \emph{Shoveling} (blue color in plots), and explicitly credit the players for finding and exploring these two different strategies. 
The difference in the two strategies as explained in \cite{sorensen2016}, is whether the tweezer is moved towards, but staying short of the initial position of the
atom, and then back towards the start and target position (yellow tunnel clan), or whether the tweezer is moved past the initial position of the atom before it is moved back towards the target position (the blue shovel clan) (Figure 3b in \cite{sorensen2016}). We note that the atom is initially caught at position 0.6.
In \cite{sorensen2016}, the two solution strategies were found by clustering the player solutions \textbf{after}
the player solutions had been improved with the KASS algorithm (the version that mitigated the sign error by not updating the amplitude as explained above).

We have shown a heat map of all the raw player solutions in Figure \ref{fig:player_heat},
and in Figure \ref{fig:yellow_blue_compare_mean_std} we have shown the mean and standard deviation of the raw player solutions that belong
to the blue and the yellow solution strategy, as tagged in \cite{sorensen2016}, respectively.
From these two figures it is clear that \textbf{the players have not explored the blue (Shoveling) strategy.}
The only strategy explored by the players, is the yellow strategy of moving up short of the atom and back again, and the main difference between the blue and the yellow groups of player solutions
is how much the amplitude is turned up and how fast that happens.
This raises a natural question: \textit{why does the blue strategy then appear when the players, in particular the player solutions from the blue group, did not explore the blue strategy at all?}
As noted above, the clustering is made on the basis of the player solutions \textbf{after} they have been processed by the KASS algorithm, and the solutions are considered at the relatively short duration of 0.17.
The KASS algorithm sweep, that is used to improve each player solution, is started from the duration considered by the player. This means that each player solution may be optimized a very different amount by the KASS algorithm
(still the version that does not change the amplitude), before they are considered by the clustering algorithm i.e. each player solution is transformed very differently (and non-linearly) before the clustering.
Concluding anything about what the players actually explore on the basis of this is misleading, in particular when a visualization of the input is not available for comparison.

\begin{figure}[ht]
  \begin{subfigure}[t]{0.48\textwidth}
    \centering
    \includegraphics[width=0.9\linewidth]{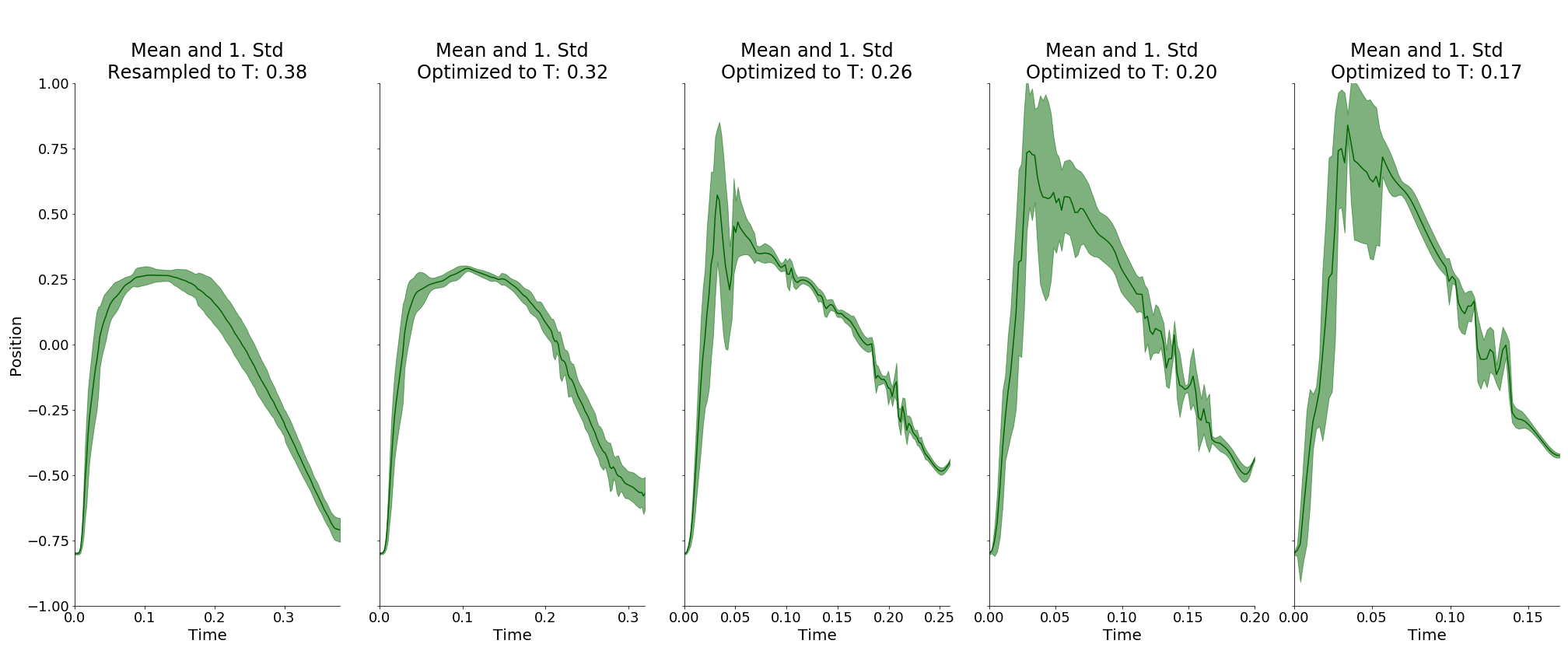}
    \caption{The development in the positions during a KASS sweep, KASS algorithm that does not change amplitudes.
      From the raw player solutions to the left, to the optimized solutions for different durations until duration 0.17 to the right. The lower the duration the more steps of the Krotov algorithm has been applied.
      Data from 10 random player solutions from the blue group.
    }
    \label{fig:yellow_blue_broken}
  \end{subfigure}%
  \hfill
  \begin{subfigure}[t]{0.48\textwidth}
    \centering
    \includegraphics[width=0.9\linewidth]{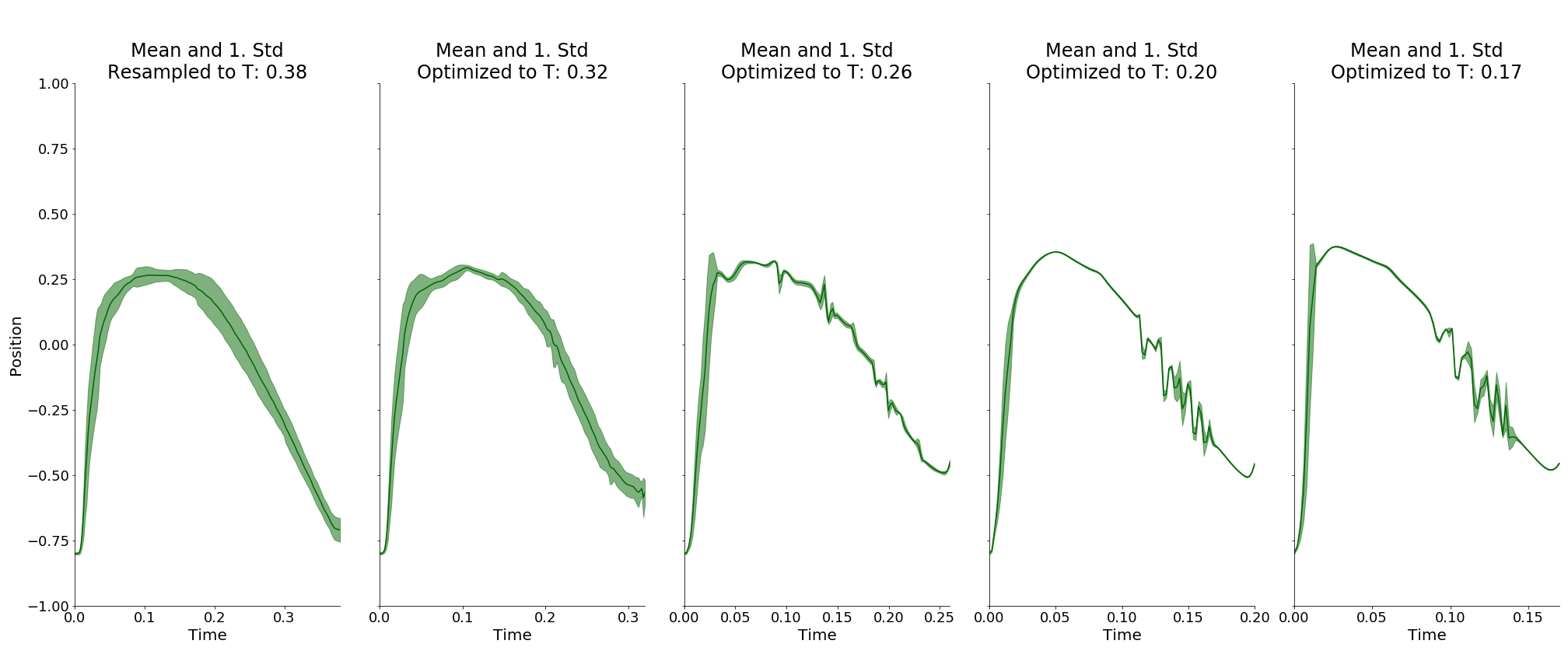}
    \caption{The development in the positions during a KASS sweep, correct KASS algorithm implementation.
      From the raw player solutions to the left, to the optimized solutions for different durations until duration 0.17 to the right. The lower the duration the more steps of the Krotov algorithm has been applied.
      Data from 10 random player solution from the blue group (same data as in Figure \ref{fig:yellow_blue_broken}. 
    }
    \label{fig:yellow_blue_fixed}
  \end{subfigure}%
  \caption{The effect of applying the KASS algorithms to the player solutions.}
\end{figure}

Now, since the players do not explore the blue strategy, it can only be the algorithm that turns the yellow player solutions into the blue solution strategies, and this is exactly what happens!
We have visualized this in Figure \ref{fig:yellow_blue_broken}. Here we can see that the KASS algorithm (the version from \cite{sorensen2016} that does not optimize amplitudes)
used on a small subset of the player solutions tagged as blue in \cite{sorensen2016},
slowly transforms the initial player solution that starts out like a yellow strategy, into a solution of the blue strategy as the algorithm iterates.
In the end when the duration is small enough, the result is  a solution following the blue strategy.
This confirms and underlines that it is a faulty KASS algorithm (it only optimizes positions) that discover and explore the blue strategy, not the players! 

This begs yet another question: \textit{why does the KASS algorithm change the player solutions so much, and why not all of them?}
The answer to this question lies mainly in the amplitudes, and how they are treated by the KASS algorithm.
As we can see from Figure \ref{fig:yellow_blue_compare_mean_std}, the main difference between the two solution strategies (the yellow and blue group player solutions),
is that the starting positions in the blue group have lower amplitude than the starting positions in the yellow.
Now, as the KASS algorithm applied to player solutions in \cite{sorensen2016} cannot change the amplitude to make up for this difference,
the KASS algorithm must find a way to change the positions to cope with the lower amplitude to get high fidelity.
In Figure \ref{fig:yellow_blue_fixed}, we have shown what happens to the exact same subset of player solutions when they are optimized with the corrected KASS algorithm
(only corrected the sign error and otherwise not changed anything from the original KASS implementation).
As the figure shows, in this case, the player solutions tagged as blue are not transformed from the yellow strategy to the blue, instead the algorithm updates the amplitudes towards -150
to match the other yellow player solutions (not shown).
This demonstrates that the purely numerical algorithm is capable of cleverly compensating in spite of the fact that it was disadvantaged from the start due to an error.
Technically speaking, the difference is that the directions followed by the KASS algorithm as dictated by the derivatives considered in the algorithm, simply changes back towards the much closer yellow local maxima when the amplitudes are updated properly.
After testing more of the player solutions tagged as blue in \cite{sorensen2016}, we found a small set of player solutions (less than 10 percent of the 100+ player solutions tagged as blue we tested)
where the corrected KASS algorithm still change the positions from  a yellow strategy to a blue. 
To investigate this further, we made an additional experiment where we increased the learning rate for the amplitude from 0.1 to 1 for the corrected KASS algorithm.
We then tested the otherwise exact same code on the same small subset of player solutions that the corrected KASS algorithm transformed from yellow to blue.
Now only one of the player solutions is changed from the yellow strategy to the blue. The player solution that the KASS algorithm with increased amplitude learning rate transforms is,
perhaps unsurprisingly, the player solution that moves closest to the atom (but still well short of).

To summarize, the results of the clustering performed in \cite{sorensen2016} is showing how the KASS algorithm behaves when it is only allowed to optimize half the parameters on very similar starting positions but with different input amplitudes,
and \textbf{not}, as concluded in \cite{sorensen2016}, that the players explore two different position strategies for the problem. They clearly do not! Furthermore, the large amount of blue player solutions found in \cite{sorensen2016} by the KASS algorithm is an artifact of the mitigation of the error of implementation used for the player solutions.
It is also clear that the clustering is highly sensitive to the exact details of the algorithm applied on top of the player solutions.

The KASS algorithm is attracted towards finding the blue strategy, but it cannot transform inputs that are too heavily biased towards the yellow solution strategy, which essentially all player solutions are.
This fits well with the findings in \cite{gronlund, Dries} that demonstrate that the standard Krotov (Figure 10 \cite{gronlund}, GRAPE (Figure 6 and Table 3 \cite{gronlund}) and Stochastic Ascent (Figure 4 and Table 1 \cite{gronlund}) algorithms, all starting from uniform random starting positions and amplitudes,
finds the blue solution strategy at the higher durations tested ($T \geq 0.16$) and uses this strategy for getting fidelity 0.999, while at lower durations the algorithms instead find the yellow solution strategy. 
Actually, around the cross-over point the algorithm may find both (Table 3 \cite{gronlund}). Furthermore, the mathematically derived protocol by D. Sels \cite{Dries, DriesarXiv},
gives a solution following the blue strategy (Figure 11 \cite{gronlund}) for duration 0.16 and it directly leads to 0.999 fidelity solutions for low durations.
The corrected KASS algorithm also finds both solution strategies, even when the algorithm has not been optimized for actually receiving correct values for the amplitude derivatives (Figure \ref{fig:kass_best_strategies} below).
Since the KASS algorithm was applied to player solutions at high durations in \cite{sorensen2016}, it is very reasonable that it changes the player solutions from the yellow to the blue strategy
if at all possible, as this is the best solution for the durations considered. 

The conclusion is that there are at least two different relevant solution strategies, they are found directly by standard algorithms for the durations where they provide the highest fidelities respectively,
and in contrast to the claims in \cite{sorensen2016}, the players only explore one of them.
The fact that player always follow the same strategy, indicates that this strategy may be inherent to the game design and how the players were taught to play it.
Whatever the reason is, the fact remains that there has been no additional insight into different strategies or exploration of the solution space by the players.
In particular, the players have not made any exploration of the blue solution strategy, which is the superior strategy when the goal is achieving 0.999 in fidelity, the stated goal of \cite{sorensen2016}. The fact that the KASS algorithm is able to transform yellow solution strategies to blue solution strategies is a testament to the algorithm, and the problem being quite easy, not the abilities of the players to explore the space of good solutions.
We note that there is one player solution that actually looks like the blue strategy in terms of positions but it has very low amplitudes.
For this reason, this player solution has essentially zero fidelity, was not tagged by the clustering as part of the blue group (nor the yellow) in \cite{sorensen2016},
and cannot be considered as an exploration of the blue strategy in any way.

\section{Learning from algorithms}
\label{sec:alg_intuition}
\begin{figure}[h]
  \begin{subfigure}[t]{0.48\textwidth}
    \centering
    \includegraphics[width=0.9\linewidth]{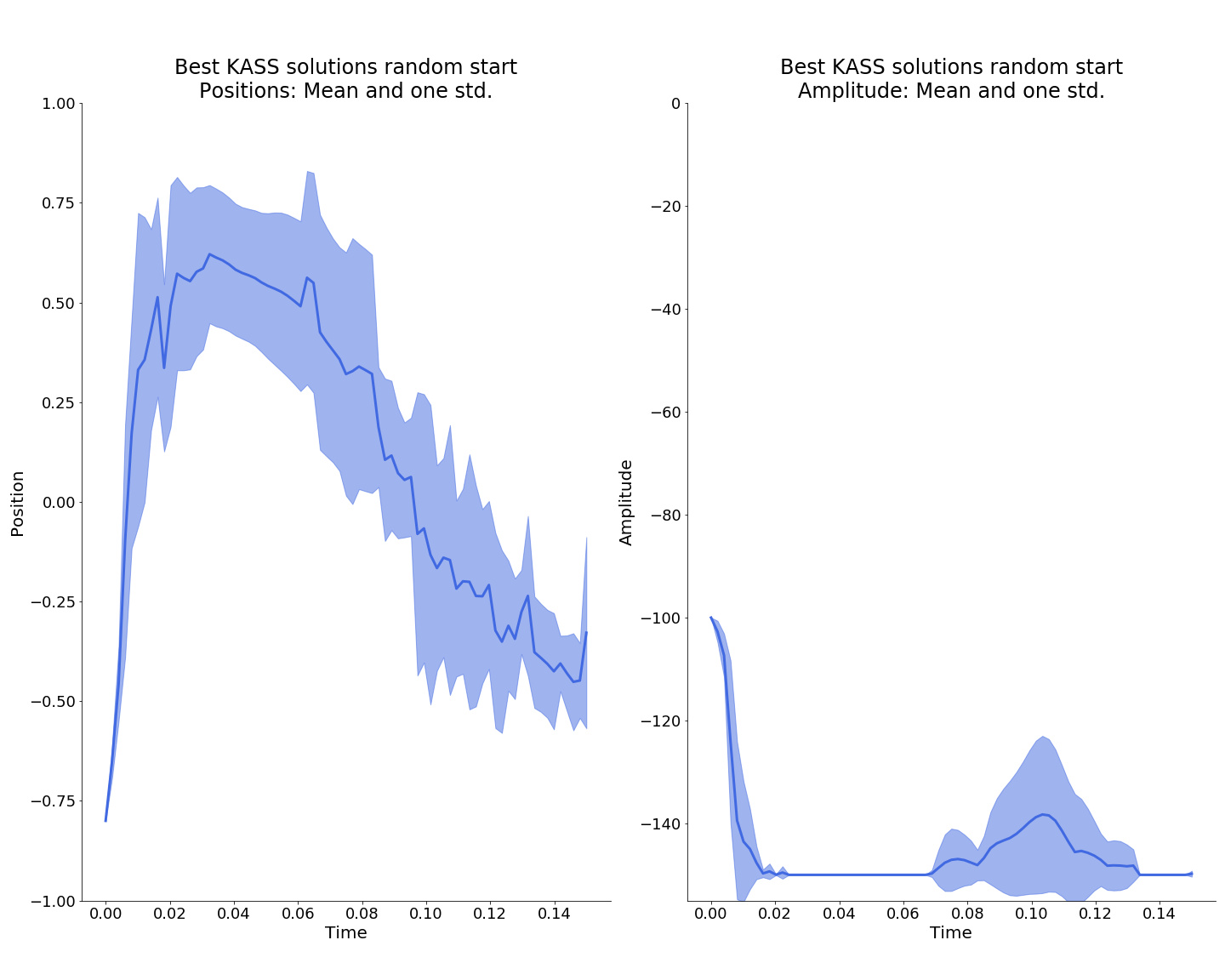}  
    \caption{The mean and 1 std. deviation of positions and amplitudes of the highest fidelity solutions found by the corrected KASS algorithm on the small subset of random starting positions and amplitudes tested.}
    \label{fig:krotov_intuition}
  \end{subfigure}%
  \hfill
  \begin{subfigure}[t]{0.48\textwidth}
    \centering
    \includegraphics[width=0.9\linewidth]{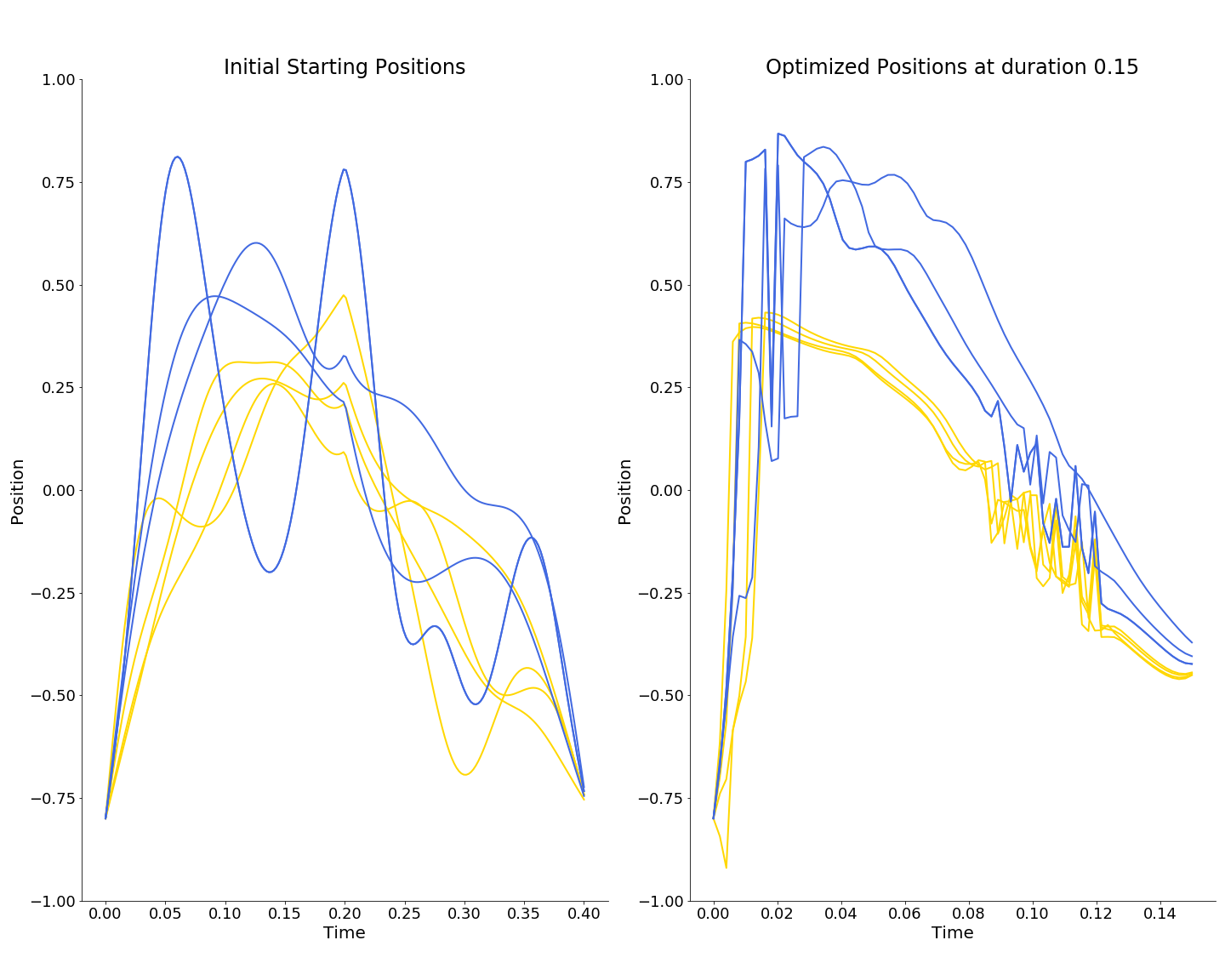}  
    \caption{The solutions (positions) that match the blue and yellow strategy found by the corrected KASS algorithm on the small subset of random starting positions and amplitudes tested.}      
    \label{fig:kass_best_strategies}
  \end{subfigure}%
  \caption{Visualization of the output for the corrected KASS algorithm}
\end{figure}

In this section we give a very short analysis of the high fidelity positions and amplitudes found by the corrected KASS algorithm.
We note that we are still using the sub-optimal hyper-parameters that were selected for the erroneous version of the KASS algorithm.
In Figure \ref{fig:krotov_intuition} we have visualized the mean position and amplitudes for each step of the best solutions found, and included the area within one standard deviation.
It is clear from this figure that the corrected KASS algorithm easily finds the solution that moves the tweezer to the atom fast and then slowly back again while shaking it, and turns the amplitude to -150. 
In Figure \ref{fig:kass_best_strategies} we have shown the subset of the best solutions found by the algorithm that match the two different solution strategies (yellow and blue).
Here we can see the best optimized solutions found among the few tested random starting positions and amplitudes include the two different strategies (yellow and blue). 
This matches the results for the other algorithm tested for BringHomeWater in \cite{gronlund}, and shows that standard algorithms automatically find the intuition behind the HILO algorithm from \cite{sorensen2016}.
This means that the players do not contribute anything to the BringHomeWater problem that is not immediately clear from applying basic algorithms.

\section{Conclusion}
In this note we have provided an in-depth analysis of the code and experiments performed in \cite{sorensen2016}.
We have shown that the only algorithm considered in \cite{sorensen2016}, the basis of the much more general conclusions made,
is incorrectly implemented, and that the conclusions from \cite{sorensen2016} are incorrect, even when considering only the non-optimized KASS algorithm in \cite{sorensen2016}.
The analysis confirms and underlines the conclusions in \cite{gronlund, Dries} that basic algorithms easily outperform all human players, naturally find the two relevant solution strategies, 
and allow for simple and highly efficient solutions for the Bring Home Water problem as considered in \cite{sorensen2016}.
Finally, the set of player solutions are highly similar variations of just one solution strategy, a solution strategy that does not yield the best QSL - the stated goal of \cite{sorensen2016},
and the players have not explored or had any insights into different strategies as otherwise claimed in \cite{sorensen2016}.

\clearpage
\appendix 

\section{Code Error}
The error in the code can be found in the file BringHomeWaterConstructor.m. In particular in line 42 where the  derivative of the Hamilton in regards to position and amplitude is implemented. 
The code is as follows:
\begin{verbatim}
bringhomewater.divpotential = {
@(u) (u(2)*4/bringhomewater.waist^2.*(bringhomewater.x - u(1))
      .*exp(-2*(bringhomewater.x - u(1)).^2/bringhomewater.waist^2)), 
@(u) (-exp(-2*(bringhomewater.x-u(1)).^2/bringhomewater.waist^2) ) 
};
\end{verbatim}
The line has two parts corresponding to the two derivatives (separated by a comma). First, we simplify the notation to work with it.

\begin{itemize}
\item $u(1)$ is the current position, we shorten it to $p$
\item $u(2)$ is the current amplitude, we shorten it to $A$
\item bringhomewater.waist = 1/4 = 0.25 a constant
\item bringhomewater.x, x is the discretization so it is a vector of 128 numbers in [-1, 1].
The code is vectorized, meaning that same computation is performed for all entries in the vector $x$ at once, so we simplify to thinking of x being a single number in [-1, 1] and it is a constant
\end{itemize}
After renaming parts,  the function we need to compute derivatives of after position $p$, and amplitude $A$ of is
$$
f(p, A) = A \exp{-32(x - p)^2}
$$
Now the derivative of $f$ after amplitude $A$ is the simplest. It is just $\exp{-32(x - p)^2}$. 
However, in the Matlab code you can see this fits perfectly except there is a minus sign in front of it!

The derivative after position $p$ is more involved, but it is
$$
(A) (\exp{- 32(x - p)^2}) (-64(x - p) (-1)) = A (4/(1/4)^2)  (x - p) (\exp{- 32(x - p)^2})
$$
exactly matching the first part of the Matlab code.
Hence, there is a sign error in the code for the two derivatives.


\begin{thebibliography}{12}
\bibitem{sorensen2016} Exploring the quantum speed limit with computer games, J.~J.~W.~H. S{\o}rensen {\it et al.},  Nature {\bf 532}, 210 (2016).
\bibitem{Dries} Stochastic gradient ascent outperforms gamers in the Quantum Moves game, D. Sels, Phys. Rev. A {\bf 97}, 040302(R) (2018).
\bibitem{krotov} Global Methods in Optimal Control Theory, V.~F. Krotov,  (Marcel Dekker, New York, 1996).
\bibitem{grape} Optimal control of coupled spin dynamics: design of nmr pulse sequences by gradient ascent algorithms, Khaneja, Reiss, Kehlet, Schulte-Herbr{\"u}ggen, Glaser, J. Magn. Reson. {\bf 172}, 296-305 (2005)
\bibitem{DriesarXiv} Can math beat gamers in Quantum Moves?, D. Sels, \url{https://arxiv.org/abs/1709.08766}
\bibitem{kingma2014adam}  Adam: A Method for Stochastic Optimization, Diederik P. Kingma and Jimmy Ba, International Conference on Learning Representations, 2015, \url{https://arxiv.org/abs/1412.6980}
\bibitem{gronlund}  Algorithms Clearly Beat Gamers at Quantum Moves. {A} Verification, Allan Gronlund, \url{http://arxiv.org/abs/1904.01008} 2019
\bibitem{molmer} (Of course) algorithms beat gamers at solving quantum control problems, Klaus Mølmer and Jacob Sherson, \url{https://arxiv.org/abs/1910.14364} 2019
  
\end{thebibliography}
\end{document}